\documentclass[12pt]{article}
\usepackage{amscd,amsfonts,amssymb,amsmath,latexsym,array,hhline}
\usepackage[dvips]{graphics}
\makeatletter
\@addtoreset{equation}{section}
\renewcommand{\@begintheorem}[2]{\begin{trivlist}\it
\item[\hspace{\labelsep}{\bf #1\ #2.}]}
\renewcommand{\@opargbegintheorem}[3]{\begin{trivlist}\it
\item[\hspace{\labelsep}{\bf #1\ #2\ (#3).}]}
\renewcommand{\@endtheorem}{\end{trivlist}}
\renewcommand{\@cite}[2]{[{#1\if@tempswa ; #2\fi}]}
\makeatother

\newcommand{\paragr}{\hspace{6mm}}
\newtheorem{Theorem}{\paragr Theorem}[section]
\newtheorem{Lemma}[Theorem]{\paragr Lemma}

\newtheorem{Proposition}[Theorem]{\paragr Proposition}
\newtheorem{Definition}[Theorem]{\paragr Definition}

\newtheorem{Example}[Theorem]{\paragr Example}
\newtheorem{Corollary}[Theorem]{\paragr Corollary}
\newcommand{\Proof}{\texttt{Proof}. }
\newcommand{\Remark}{\texttt{Remark}. }
\newcommand{\Iremark}{\texttt{Important Remark}. }

\newcommand{\Remarks}{\texttt{Remarks}. }

\newcommand{\End}{~\hfill $\Box\!$}

\newcommand{\skm}{\medskip}
\newcommand{\skb}{\bigskip}
\newcommand{\sksminus}{}
\newcommand{\sksmin}{}

\newcommand{\al}{\alpha}

\newcommand{\ga}{\gamma}

\newcommand{\de}{\delta}

\newcommand{\la}{\lambda}
\newcommand{\si}{\sigma}

\newcommand{\eps}{\varepsilon}
\renewcommand{\phi}{\varphi}
\renewcommand{\kappa}{\varkappa}
\newcommand{\N}{\mathbb{N}}

\newcommand{\R}{\mathbb{R}}

\newcommand{\F}{\mathcal{F}}

\newcommand{\DDD}{\mathcal{D}}

\newcommand{\MMM}{\mathcal{M}}
\newcommand{\RRR}{\mathcal{R}}
\newcommand{\PPP}{\mathcal{P}}
\newcommand{\EX}{\mathcal{X}}
\newcommand{\PD}{\mathcal{PD}}
\newcommand{\RD}{\mathcal{RD}}
\newcommand{\EE}{\mathsf{E}}
\newcommand{\PP}{\mathsf{P}}
\newcommand{\QQ}{\mathsf{Q}}

\newcommand{\emp}{\emptyset}
\newcommand{\lb}{\langle}
\newcommand{\rb}{\rangle}
\newcommand{\wt}{\widetilde}
\newcommand{\wl}{\overline}
\newcommand{\xra}{\xrightarrow}
\newcommand{\da}{\downarrow}

\newcommand{\Lea}{\Longleftrightarrow}
\newcommand{\ds}{\displaystyle}

\newcommand{\Law}{\mathop{\rm Law}\nolimits}
\newcommand{\conv}{\mathop{\rm conv}\nolimits}

\newcommand{\supp}{\mathop{\rm supp}\nolimits}
\newcommand{\pr}{\mathop{\rm pr}\nolimits}

\newcommand{\esssup}{\mathop{\rm esssup}}
\newcommand{\essinf}{\mathop{\rm essinf}}
\newcommand{\argmax}{\mathop{\rm argmax}}
\newcommand{\argmin}{\mathop{\rm argmin}}

\newcommand{\cl}{\mathop{\rm cl}}
\newcommand{\cov}{\mathop{\sf cov}}
\newcommand{\var}{\mathop{\sf var}}
\newcommand{\corr}{\mathop{\sf corr}}
\renewcommand{\inf}{\mathop{\rm inf\rule[-0.8mm]{0mm}{1mm}}}

\newcommand{\RAROC}{{\rm RAROC}}
\newcommand{\VaR}{{\rm V@R}}

\newcommand{\INGD}{I_{\text{\sl NGD}}}
\newcommand{\INGDR}{I_{\text{\sl NGD}\,(\!R)}}
\newcommand{\INA}{I_{\text{\sl NA}}}

\newenvironment{mitemize}%
{\begin{list}{$\bullet$}{
\leftmargin=32pt
\rightmargin=0pt
\labelsep=5pt
\labelwidth=20pt
\itemindent=0pt
\topsep=5pt plus 2pt minus 4pt
\partopsep=2pt plus 1pt minus 1pt
\parsep=0pt
\itemsep=0pt}}%
{\end{list}}

\begin{document}
\vspace*{1mm}
\begin{center}\bf
PRICING WITH COHERENT RISK
\end{center}

\begin{center}\itshape\bfseries
A.S.~Cherny
\end{center}

\begin{center}
\textit{Moscow State University,}\\
\textit{Faculty of Mechanics and Mathematics,}\\
\textit{Department of Probability Theory,}\\
\textit{119992 Moscow, Russia.}\\
\texttt{E-mail: cherny@mech.math.msu.su}\\
\texttt{Webpage: http://mech.math.msu.su/\~{}cherny}
\end{center}

\begin{abstract}
\textbf{Abstract.}
This paper deals with applications of coherent
risk measures to pricing in incomplete markets. Namely, we study
the No Good Deals pricing technique based on coherent risk. Two
forms of this technique are presented: one defines a good deal as
a trade with negative risk; the other one defines a good deal as a
trade with unusually high RAROC. For each technique, the
fundamental theorem of asset pricing and the form of the
fair price interval are presented. The model considered
includes static as well as dynamic models, models with an
infinite number of assets, models with transaction costs,
and models with portfolio constraints.
In particular, we prove that in a model with proportional
transaction costs the fair price interval converges to the
fair price interval in a frictionless model as the
coefficient of transaction costs tends to zero.

Moreover, we study some problems in the ``pure'' theory
of risk measures: we present a simple geometric solution
of the capital allocation problem and apply it to define
the coherent risk contribution.

The mathematical tools employed are probability theory,
functional analysis, and finite-dimensional convex
analysis.

\bigskip
\textbf{Key words and phrases:}
Capital allocation,
coherent risk measures,
extreme measures,
generating set,
No Good Deals,
RAROC,
risk contribution,
risk-neutral measures,
support function,
Tail V@R,
transaction costs,
Weighted V@R.
\end{abstract}

\section{Introduction}
\label{I}

\textbf{1. Overview.}
The three basic pillars of finance are:
\begin{mitemize}
\item optimal investment;
\item pricing and hedging;
\item risk measurement and management.
\end{mitemize}
The most well-known financial theories related to the
first pillar are the Markowitz mean-variance analysis
and Sharpe's CAPM, which are often termed the
``first revolution in finance''.
The most well-known result related to the second pillar
is the Black--Scholes--Merton formula, which is often
termed the ``second revolution in finance''.
Recently a very important innovation has appeared in
connection with the third pillar.
In 1997, Artzner, Delbaen, Eber and Heath~\cite{ADEH97},
\cite{ADEH99} introduced the concept of a
\textit{coherent risk measure} as a new way of measuring risk.
Since 1997, the theory of coherent risk measures
has rapidly been evolving and is already termed in some sources
the ``third revolution in finance'' (see~\cite{Sz04}).
Let us mention, in particular, the papers~\cite{A02},
\cite{AT02}, \cite{D02}, \cite{FS02}, \cite{FS02b},
\cite{JMT04}, \cite{K01}, \cite{T02}
and the reviews~\cite{D05}, \cite[Ch.~4]{FS04}, \cite{S05}.
Currently, one of the major tasks is the problem
of proper risk measurement in the dynamic setting;
see, in particular,~\cite{CDK05}, \cite{DS05},
\cite{JR05}, \cite{R04}, and~\cite{RSE05}.

The theory of coherent risk measures is important not only
for risk measurement.
Indeed, risk ($\approx$ uncertainty) is at the very basis
of the whole finance, and therefore, a new way of looking
at risk yields new approaches to other problems of
finance, in particular, to those related to the first
and the second pillars.
Nowadays, more and more research is aimed at
\textit{applications} of coherent risk measures to other
problems of finance.

One of the major goals of modern financial mathematics
is providing adequate price bounds for
derivative contracts in incomplete markets.
It is known that No Arbitrage price bounds in incomplete
markets are typically unacceptably wide, and fundamentally
new ideas are required to narrow these bounds.
Recently, a promising approach to this problem termed
\textit{No Good Deals} (\textit{NGD}) pricing has been
proposed in~\cite{BL00}, \cite{CS00}.
Let us illustrate its idea by an example. Consider a
contract that with probability $1/2$ yields nothing
and with probability $1/2$ yields 1000~USD.
The No Arbitrage (NA) price interval for this contract
is $(0,1000)$. But if the price of the contract is, for
instance, 15~USD, then everyone would be willing to buy it,
and the demand would not match the supply.
Thus, 15~USD is an unrealistic price because it yields a
good deal, i.e. a trade that is attractive to most market
participants. The technique of the NGD pricing
is based on the assumption that good deals do not exist.

A problem that arises immediately is how to define
a good deal. There is no canonical answer, and several
approaches have been proposed in the literature.
Cochrane and Sa\'a-Requejo~\cite{CS00}
defined a good deal as a trade with unusually high
Sharpe ratio,
Bernardo and Ledoit~\cite{BL00} based their definition
on another gain to loss ratio,
while \v{C}ern\'y and Hodges~\cite{CH01}
proposed a generalization of both definitions
(see also the paper~\cite{BS05} by Bjork and Slinko,
which extends the results of~\cite{CS00}).

The technique of the NGD pricing can also be motivated as
follows. When a trader sells a contract, he/she would charge for it
a price, with which he/she will be able to superreplicate the contract.
In theory the superreplication is typically understood
almost surely, but in practice an agent looks for an offsetting
position such that the risk of his/her overall portfolio would
stay within the limits prescribed by his/her management
(the almost sure superreplication is virtually
impossible in practice).
These considerations lead to the NGD pricing with a good
deal defined as a trade with negative risk.
Now, if risk is measured by \VaR, this technique leads to
the quantile hedging introduced by F\"ollmer and Leukert~\cite{FL99}.
But instead of \VaR, one can take a coherent risk measure.
The corresponding pricing technique has already been
considered in several papers.
Carr, Geman, and Madan~\cite{CGM01}
(see also the review paper~\cite{CGM02})
studied this technique in a probabilistic framework
(although they do not use the term ``good deal''), while
Jaschke and K\"uchler~\cite{JK01}
studied this technique in a topological space framework
in the spirit of Harrison and Kreps~\cite{HK79}
(see also the paper~\cite{S04} by Staum, which extends the
results of~\cite{JK01}).
Furthermore, Larsen, Pirvu, Shreve, and
T\"ut\"unc\"u~\cite{LPST04} considered pricing based on
convex risk measures instead of coherent ones
(convex risk measures were introduced by
F\"ollmer and Schied~\cite{FS02}).
Roorda, Schumacher, and Engwerda~\cite{RSE05}
studied pricing in the multiperiod model using as a
basis dynamic coherent risk measures instead of static
ones.

\skb
\textbf{2. Goal of the paper.}
This is the first of a series of papers dealing with
applications of coherent risk measures to the basic
problems of finance (the other paper in the series
is~\cite{C062}). The basic idea behind the series is:
\begin{center}\bfseries\itshape
the whole finance can be built based on coherent risks.
\end{center}
In this paper, we study
applications to pricing in incomplete markets.
Our approach is similar to that of~\cite{CGM01},
but~\cite{CGM01} assumes an unrealistic world of a finite
state space and a finite set of probabilistic scenarios
defining a coherent risk measure (most natural coherent
risk measures are defined through an infinite set of
probabilistic scenarios; see Subsection~\ref{BD}).
Our model is general in the sense that we consider
an arbitrary~$\Omega$ and a general class of coherent
risk measures (satisfying only a sort of compactness condition).
Moreover, our approach applies to dynamic models,
to models with an infinite number of assets,
to models with transaction costs,
and to models with convex portfolio constraints.
Within this general model, we prove the Fundamental
Theorem of Asset Pricing (Theorem~\ref{UP4}) and provide
the form of the fair price interval of a contingent claim
(Corollary~\ref{UP6}).
We confine ourselves to static risk measures.

A problem that has attracted attention in several papers
is as follows. Consider a model with proportional
transaction costs. Is it true that the upper (resp.,
lower) price of a contingent claim in this model
tends to the upper (resp., lower) price of this claim
in the frictionless model as the coefficient of transaction
costs tends to zero?
It was shown in~\cite{C05b}, \cite{CPT99}, \cite{LS97},
and~\cite{SSC95} that, for NA prices, the answer to this
question is negative already in the Black-Scholes model
(the contingent claim considered in these papers is a
European call option).
This result might be interpreted as follows:
the NA technique is useless in continuous-time models
with transaction costs.
In this paper (Theorem~\ref{PDT2}),
we prove that, for NGD prices, the answer
to the above question is positive.
This is done within a framework of a general model
(the price follows an arbitrary process)
with an infinite number of assets and an arbitrary
contingent claim (satisfying only some integrability
condition).
The advantage of the NGD pricing is not only that this
result is true, but also that its proof is very short.

Furthermore, we introduce a new variant of pricing
based on coherent risk, which we call the
\textit{RAROC-based NGD} pricing. The idea is to define
a good deal as a trade with unusually high
Risk-Adjusted Return on Capital (RAROC), where RAROC
is defined through coherent risk.
On the mathematical side, this technique is reduced
to the standard NGD pricing (with the original risk
measure replaced by another one).

Although this series of papers deals primarily with
applications of coherent risk measures to problems of
finance, we also establish some results and give
several definitions related to ``pure'' risk measures
(these are needed for applications).
In particular, we introduce the notion of an
\textit{extreme measure}. The results of this paper
and~\cite{C062} show that this notion is very convenient
and important; it appears in the outcomes of several
pricing techniques proposed in~\cite{C062} and in
considerations of the equilibrium problem in~\cite{C062}.
In the present paper, we provide a solution of the
\textit{capital allocation} problem in terms of extreme
measures (Theorem~\ref{CA1}).
Let us remark that this problem was considered
in~\cite{D05}, \cite{D01}, \cite{F03}, \cite{K05},
\cite{O99}, and~\cite{T02}.

Parallel with the measurement of outstanding risks, a
very important problem is measuring
the risk contribution of a subportfolio to a ``big''
portfolio. Based on our solution of the capital allocation
problem, we propose several equivalent definitions of the
coherent \textit{risk contribution}.

Another notion we introduce is the notion of a
\textit{generator}. It establishes a bridge between
coherent risks and convex analysis, opening the way for
geometry. In particular, we provide (see Figure~1) a
geometric solution of the capital allocation problem
(thus there are two solutions: a probabilistic one is
given in terms of extreme measures, while a geometric one
is given in terms of generators).
We also provide a geometric solution of the
pricing and hedging problem (Proposition~\ref{HI4}) for a
model with a finite number of assets.
Furthermore, we provide in~\cite{C062} geometric solutions
of several optimization problems, optimality pricing
problems, and the equilibrium problem.
In fact, for most problems considered in this series
of papers, we provide two sorts of results:
\begin{mitemize}
\item a geometric result applicable to a model with a finite
number of assets is given in terms of generators;
\item a probabilistic result applicable to a general model
is typically given in terms of extreme measures.
\end{mitemize}

\skm
\textbf{3. Structure of the paper.}
Section~\ref{M} deals with ``pure'' risk measures rather
than with their applications.
Subsection~\ref{BD} recalls some basic definitions
related to coherent risks.
In Subsection~\ref{LO}, we introduce the $L^1$-spaces
associated with a coherent risk measure (these are employed
in the technical conditions in theorems below).
Subsection~\ref{EM} presents the definition of an
extreme measure.
In Subsection~\ref{CA}, we provide a solution of the
capital allocation problem.
Subsection~\ref{RC} deals with equivalent definitions
of risk contribution.

Section~\ref{P} is related to the NGD pricing.
In Subsections~\ref{UP} and~\ref{RP}, we study the
ordinary and the RAROC-based forms of this technique,
respectively. The model considered is a general one, and in
Subsections~\ref{PF}--\ref{PDT}, we consider some
particular cases of this model:
a static model with a finite number of assets
(for which fair price intervals admit a simple geometric
description; see Figure~3),
a continuous-time dynamic model, and
a continuous-time dynamic model with transaction costs.
Furthermore, in Subsection~\ref{HI}, we provide a
geometric solution of the hedging problem for a static
model with a finite number of assets (see Figure~4).

\bigskip
{\itshape\bfseries Acknowledgement.}
I am thankful to D.B.~Madan for valuable discussions
and important advice.

\section{Coherent Risk Measures}
\label{M}

\subsection{Basic Definitions}
\label{BD}

Let $(\Omega,\F,\PP)$ be a probability space.
The following definition was introduced in~\cite{ADEH97},
\cite{ADEH99}. These papers considered only a finite~$\Omega$,
in which case the continuity axiom~(e) is not needed.
It was added for a general~$\Omega$ by Delbaen~\cite{D02}.

\begin{Definition}\rm
\label{BD1}
A \textit{coherent utility function on $L^\infty$} is a
map $u:L^\infty\to\R$ with the properties:
\begin{mitemize}
\item[(a)] (Superadditivity) $u(X+Y)\ge u(X)+u(Y)$;
\item[(b)] (Monotonicity) If $X\le Y$, then $u(X)\le u(Y)$;
\item[(c)] (Positive homogeneity) $u(\la X)=\la u(X)$ for $\la\in\R_+$;
\item[(d)] (Translation invariance) $u(X+m)=u(X)+m$ for $m\in\R$;
\item[(e)] (Fatou property) If $|X_n|\le1$,
$X_n\xra{\PP}X$, then $u(X)\ge\limsup_n u(X_n)$.
\end{mitemize}
The corresponding \textit{coherent risk measure} is
$\rho(X)=-u(X)$.
\end{Definition}

\Remark
Typically, a coherent risk measure is defined only via
conditions (a)--(d), and then one speaks about coherent
risk measures with the Fatou property.
However, only such risk measures are useful, and for this
reason we find it more convenient to add~(e) as a basic axiom.

\skm
The theorem below was established in~\cite{ADEH99} for
the case of a finite $\Omega$ (in this case the axiom~(e)
is not needed) and in~\cite{D02} for the general case.
We denote by $\PPP$ the set of probability measures
on~$\F$ that are absolutely continuous with respect to~$\PP$.
Throughout the paper, we identify measures from~$\PPP$
(these are typically denoted by~$\QQ$)
with their densities with respect to~$\PP$
(these are typically denoted by~$Z$).

\begin{Theorem}[Basic representation theorem]
\label{BD2}
A function~$u$ satisfies
conditions {\rm(a)--(e)} if and only if
there exists a nonempty set $\DDD\subseteq\PPP$ such that
\begin{equation}
\label{bd1}
u(X)=\inf_{\QQ\in\DDD}\EE_\QQ X,\quad X\in L^\infty.
\end{equation}
\end{Theorem}

So far, a coherent risk measure has been defined on bounded
random variables.
Let us ask ourselves the following question:
Are ``financial'' random variables like the increment
of a price of some asset indeed bounded?
The right way to address this question is to
split it into two parts:
\begin{mitemize}
\item Are ``financial'' random variables bounded in practice?
\item Are ``financial'' random variables bounded in theory?
\end{mitemize}
The answer to the first question is positive
(clearly, everything is bounded by the number of the
atoms in the universe).
The answer to the second question is negative because
most distributions used in theory (like the lognormal
one) are unbounded.
So, as we are dealing with theory, we need to extend
coherent risk measures to the space $L^0$ of all random
variables.
It is hopeless to axiomatize the notion of a risk measure
on $L^0$ and then to obtain the corresponding
representation theorem.
Instead, we take representation~\eqref{bd1} as the
basis and extend it to~$L^0$.

\begin{Definition}\rm
\label{BD3}
A \textit{coherent utility function on $L^0$} is a map
$u:L^0\to[-\infty,\infty]$ defined~as
\begin{equation}
\label{bd2}
u(X)=\inf_{\QQ\in\DDD}\EE_\QQ X,\quad X\in L^0,
\end{equation}
where $\DDD\subseteq\PPP$ and $\EE_\QQ X$ is understood as
$\EE_\QQ X^+-\EE_\QQ X^-$ with the convention
$\infty-\infty=-\infty$.
The corresponding \textit{coherent risk measure} is
$\rho(X)=-u(X)$.
\end{Definition}

Clearly, a set $\DDD$, for which representations~\eqref{bd1}
and~\eqref{bd2} are true, is not unique.
However, there exists the largest such set given by
$\{\QQ\in\PPP:\EE_\QQ X\ge u(X)\text{ for any }X\}$.
We introduce the following definition.

\begin{Definition}\rm
\label{BD4}
We will call the largest set, for which~\eqref{bd1}
{\rm(}resp., \eqref{bd2}{\rm)} is true, the
\textit{determining set of}~$u$.
\end{Definition}

\Remark
Clearly, the determining set is convex.
For coherent utility functions on~$L^\infty$, it is also
$L^1$-closed.
However, for coherent utility functions on $L^0$, it is
not necessarily $L^1$-closed.
As an example, take a positive unbounded random
variable $X_0$ such that $\PP(X_0=0)>0$ and consider
$\DDD_0=\{\QQ\in\PPP:\EE_\QQ X_0=1\}$.
Clearly, the determining set $\DDD$ of the coherent
utility function $u(X)=\inf_{\QQ\in\DDD_0}\EE_\QQ X$
satisfies
$\DDD_0\subseteq\DDD\subseteq\{\QQ\in\PPP:\EE_\QQ X_0\ge1\}$.
On the other hand, the $L^1$-closure of $\DDD_0$ contains
a measure $\QQ_0$ concentrated on $\{X_0=0\}$.

\skm
\Iremark
Let $\DDD$ be an $L^1$-closed convex subset
of~$\PPP$.
(Let us note that a particularly important case is
where $\DDD$ is $L^1$-closed, convex, and uniformly
integrable; this condition will be needed in a number
of places below).
Define a coherent utility function~$u$
by~\eqref{bd2}. Then $\DDD$ is the determining set of~$u$.
Indeed, assume that the determining set $\wt\DDD$ is
greater than~$\DDD$, i.e. there exists
$\QQ_0\in\wt\DDD\setminus\DDD$. Then, by the Hahn-Banach
theorem, we can find $X_0\in L^\infty$ such that
$\EE_{\QQ_0}X_0<\inf_{\QQ\in\DDD}\EE_\QQ X$, which is a
contradiction.
The same argument shows that $\DDD$ is also the determining
set of the restriction of~$u$ to $L^\infty$.

\skm
In what follows, we will always consider coherent utility
functions on~$L^0$.

\begin{Example}\rm
\label{BD5}
{\bf(i)} \textit{Tail V@R} (the terms
\textit{Average V@R}, \textit{Conditional V@R},
and \textit{Expected Shortfall} are also used)
is the risk measure corresponding to the
coherent utility function
$$
u_\la(X)=\inf_{\QQ\in\DDD_\la}\EE_\QQ X,
$$
where $\la\in[0,1]$ and
\begin{equation}
\label{bd3}
\DDD_\la=\biggl\{\QQ\in\PPP:\frac{d\QQ}{d\PP}\le\la^{-1}\biggr\}.
\end{equation}
In particular, if $\la=0$, then the corresponding
coherent utility function has the form
$u(X)=\essinf_\omega X(\omega)$.
For more information on Tail V@R,
see~\cite{AT02}, \cite[Sect.~6]{D02}, \cite[Sect.~7]{D05},
\cite[Sect.~4.4]{FS04}, \cite[Sect.~1.3]{S05}.

{\bf(ii)} \textit{Weighted V@R on $L^\infty$}
(the term \textit{spectral risk measure} is also used)
is the risk measure corresponding to the coherent
utility function
$$
u_\mu(X)=\int_{[0,1]}u_\la(X)\mu(d\la),\quad X\in L^\infty,
$$
where $\mu$ is a probability measure on $[0,1]$.

\textit{Weighted V@R on $L^0$} is the risk
measure corresponding to the coherent utility function
$$
u_\mu(X)=\inf_{\QQ\in\DDD_\mu}\EE_\QQ X,\quad X\in L^0,
$$
where $\DDD_\mu$ is the determining set of $u_\mu$
on $L^\infty$.

Let us remark that, under some regularity conditions
on~$\mu$, Weighted V@R possesses some nice properties
that are not shared by Tail V@R.
In a sense, it is ``smoother'' than Tail V@R.
We consider Weighted V@R as one of the most important
classes (or maybe the most important class) of coherent
risk measures.
For a detailed study of this risk measure,
see~\cite{A02}, \cite{A04}, \cite{Dowd05},
\cite{K01} as well as the paper~\cite{C05e},
which is in some sense the continuation of the present paper.\End
\end{Example}

\subsection{Spaces $L_w^1$ and $L_s^1$}
\label{LO}

For a subset $\DDD$ of $\PPP$, we introduce the
\textit{weak} and \textit{strong $L^1$-spaces}
\begin{align*}
L_w^1(\DDD)&=\{X\in L^0:u(X)>-\infty,\,u(-X)>-\infty\},\\
L_s^1(\DDD)&=\Bigl\{X\in L^0:\lim_{n\to\infty}
\sup_{\QQ\in\DDD}\EE_\QQ|X|I(|X|>n)=0\Bigr\}.
\end{align*}
Clearly, $L_s^1(\DDD)\subseteq L_w^1(\DDD)$.
If $\DDD=\{\QQ\}$ is a singleton, then
$L_w^1(\DDD)=L_s^1(\DDD)=L^1(\QQ)$, which motivates
the notation.

In general, $L_s^1(\DDD)$ might be strictly smaller than
$L_w^1(\DDD)$.
Indeed, let $X_0$ be a positive unbounded random variable
with $\PP(X_0=0)>0$ and let $\DDD=\{\QQ\in\PPP:\EE_\QQ X_0=1\}$.
Then $X_0\in L_w^1(\DDD)$, but $X_0\notin L_s^1(\DDD)$.
(One can also construct a similar counterexample with
an $L^1$-closed set $\DDD$; see Example~\ref{EM3}).
However, as shown by the proposition below, in most natural
situations weak and strong $L^1$-spaces coincide.

\begin{Proposition}
\label{LO1}
{\bf(i)} If $\DDD_\la$ is the determining set of Tail V@R
{\rm(}see Example~\ref{BD5}~{\rm(}i{\rm))}
with $\la\in(0,1]$, then $L_w^1(\DDD_\la)=L_s^1(\DDD_\la)$.

{\bf(ii)} If $\DDD_\mu$ is the determining set of Weighted V@R
{\rm(}see Example~\ref{BD5}~{\rm(}ii{\rm))}
with $\mu$ concentrated on $(0,1]$, then
$L_w^1(\DDD_\mu)=L_s^1(\DDD_\mu)$.

{\bf(iii)} If all the densities from $\DDD$ are bounded
by a single constant and $\PP\in\DDD$,
then $L_w^1(\DDD)=L_s^1(\DDD)$.

{\bf(iv)} If $\DDD$ is a convex combination
$\sum_{n=1}^N a_n\DDD_n$, where $\DDD_1,\dots,\DDD_N$ are such
that $L_w^1(\DDD_n)=L_s^1(\DDD_n)$,
then $L_w^1(\DDD)=L_s^1(\DDD)$.

{\bf(v)} If $\DDD=\conv(\DDD_1,\dots,\DDD_N)$,
where $\DDD_1,\dots,\DDD_N$ are such
that $L_w^1(\DDD_n)=L_s^1(\DDD_n)$,
then $L_w^1(\DDD)=L_s^1(\DDD)$.
\end{Proposition}

\begin{Lemma}
\label{LO2}
If $\mu$ is a convex combination
$\sum_{n=1}^\infty a_n\de_{\la_n}$, where $\la_n\in(0,1]$,
then the determining set $\DDD_\mu$ of Weighted V@R
corresponding to~$\mu$ has the form
$\sum_{n=1}^\infty a_n\DDD_{\la_n}$,
where $\DDD_\la$ is given by~\eqref{bd3}.
\end{Lemma}

\sksminus
\Proof
Denote $\sum_n a_n\DDD_{\la_n}$ by $\DDD$.
Clearly, $\DDD$ is convex.
Fix $X\in L^\infty$.
It is easy to see that, for any~$n$, the minimum of
expectations of $\EE XZ$ over $Z\in\DDD_{\la_n}$ is
attained (for more details, see~\cite[Prop.~2.7]{C05e}).
Hence, the minimum of expectations $\EE_\PP XZ$
over $Z\in\DDD$ is attained.
By the James theorem (see~\cite{F80}), $\DDD$ is weakly
compact.
As it is convex, an application of the Hahn-Banach
theorem shows that it is $L^1$-closed.

Obviously, $u_\mu(X)=\inf_{\QQ\in\DDD}\EE_\QQ X$ for any
$X\in L^\infty$.
Taking into account the Important Remark following
Definition~\ref{BD4}, we get $\DDD_\mu=\DDD$.\End

\skm
\texttt{Proof of Proposition~\ref{LO1}.}
The only nontrivial statement is (ii).
In order to prove it, consider the measures
$\wt\mu=\sum_{k=1}^\infty a_k\de_{2^{-k}}$,
$\bar\mu=\sum_{k=1}^\infty a_k\de_{2^{-k+1}}$,
where $a_k=\mu((2^{-k},2^{-k+1}])$.
As $u_{\wt\mu}\le u_\mu\le u_{\bar\mu}$, we have
$\DDD_{\wt\mu}\supseteq\DDD_\mu\supseteq\DDD_{\bar\mu}$.
By Lemma~\ref{LO2},
$$
\DDD_{\wt\mu}=\biggl\{\sum_{k=1}^\infty a_k Z_k:
Z_k\in\DDD_{2^{-k}}\biggr\},\qquad
\DDD_{\bar\mu}=\biggl\{\sum_{k=1}^\infty a_k Z_k:
Z_k\in\DDD_{2^{-k+1}}\biggr\}.
$$
Take $X\in L_w^1(\DDD_\mu)$.
Consider $Z_k=2^{k-1}I(X<q_k)+c_k I(X=q_k)$, where $q_k$ is the
$2^{-k+1}$-quantile of $X$ and $c_k$ is chosen in such a way
that $\EE_\PP Z_k=1$. Then
$$
\EE_\PP Z_k X=\min_{Z\in\DDD_{2^{-k+1}}}\EE_\PP ZX.
$$
The density $Z_0=\sum_{k=1}^\infty a_kZ_k$ belongs to $\DDD_{\bar\mu}$ and
$$
\EE_\PP Z_0X=\min_{Z\in\DDD_{\bar\mu}}\EE_\PP ZX.
$$
In view of the inclusion
$X\in L_w^1(\DDD_\mu)\subseteq L_w^1(\DDD_{\bar\mu})$,
the latter quantity is finite. Thus,
$$
\sum_{k=1}^\infty a_k\min_{Z\in\DDD_{2^{-k+1}}}\EE_\PP ZX>-\infty,
$$
which implies that
$$
\sum_{k=1}^\infty a_k\min_{Z\in\DDD_{2^{-k+1}}}\EE_\PP Z(-X^-)>-\infty.
$$
The same estimate is true for $X^+$, and therefore,
\begin{equation}
\label{em2}
\sum_{k=1}^\infty a_k\sup_{Z\in\DDD_{2^{-k}}}\EE_\PP Z|X|
\le2\sum_{k=1}^\infty a_k\sup_{Z\in\DDD_{2^{-k+1}}}\EE_\PP Z|X|
<\infty.
\end{equation}
It is clear that $X\in L^1$, and thus, for each~$k$,
$$
\sup_{Z\in\DDD_{2^{-k}}}\EE_\PP Z|X|I(|X|>n)
\le 2^k\EE_\PP|X|I(|X|>n)\xra[n\to\infty]{}0.
$$
This, combined with~\eqref{em2}, yields
\begin{align*}
\sup_{Z\in\DDD_\mu}\EE_\PP Z|X|I(|X|>n)
&\le\sup_{Z\in\DDD_{\wt\mu}}\EE_\PP Z|X|I(|X|>n)\\
&=\sum_{k=1}^\infty a_k\sup_{Z\in\DDD_{2^{-k}}}\EE_\PP Z|X|I(|X|>n)
\xra[n\to\infty]{}0.
\end{align*}

\subsection{Extreme Measures}
\label{EM}

\begin{Definition}\rm
\label{EM1}
Let $u$ be a coherent utility function with the
determining set~$\DDD$. Let $X\in L^0$.
We will call a measure $\QQ\in\DDD$ an
\textit{extreme measure} for~$X$ if
$\EE_\QQ X=u(X)\in(-\infty,\infty)$.

The set of extreme measures will be denoted by $\EX_\DDD(X)$.
\end{Definition}

Let us recall some general facts related to the
weak topology on~$L^1$.
The \textit{weak topology} on $L^1$ is
induced by the duality between $L^1$ and $L^\infty$
and is usually denoted as $\si(L^1,L^\infty)$.
The Dunford-Pettis criterion states that a set
$\DDD\subseteq\PPP$ is weakly compact if and only if it
is weakly closed and uniformly integrable.
Furthermore, an application of the Hahn-Banach theorem
shows that a convex set $\DDD\subseteq\PPP$ is weakly
closed if and only if it is $L^1$-closed.

\begin{Proposition}
\label{EM2}
If the determining set $\DDD$ is weakly compact and $X\in L_s^1(\DDD)$,
then $\EX_\DDD(X)\ne\emp$.
\end{Proposition}

\sksminus
\Proof
It is clear that $u(X)\in(-\infty,\infty)$.
Find a sequence $Z_n\in\DDD$ such that
$\EE_\PP Z_n X\to u(X)$.
This sequence has a weak limit point $Z_\infty\in\DDD$.
Clearly, the map $\DDD\ni Z\mapsto\EE_\PP ZX$ is
weakly continuous.
Hence, $\EE_\PP Z_\infty X=u(X)$, which means that
$Z_\infty\in\EX_\DDD(X)$.\End

\begin{Example}\rm
\textbf{(i)} If $u$ corresponds to Tail V@R of order
$\la\in(0,1]$ (see Example~\ref{BD5}~(i))
and $X$ has a continuous distribution,
then it is easy to see that $\EX_\DDD(X)$ consists of a
unique density $\la^{-1}I(X\le q_\la)$, where $q_\la$
is a $\la$-quantile of~$X$.

\textbf{(ii)} If $u$ corresponds to Weighted V@R with
the weighting measure~$\mu$ (see Example~\ref{BD5}~(ii))
and $X$ has a continuous distribution,
then $\EX_\DDD(X)$ consists of a unique density
$g(X)$, where $g(x)=\int_{[F(x),1]}\la^{-1}\mu(d\la)$
and $F$ is the distribution function of~$X$
(see~\cite[Sect.~6]{C05e}).
Note that this density reflects the risk aversion of an
agent possessing a portfolio that produces the
P\&L (Profit\&Loss)~$X$.\End
\end{Example}

\skm
The condition that $\DDD$ should be weakly compact is
very mild and is satisfied for the determining sets of
most natural coherent risk measures.
For example, the determining set $\DDD_\la$ of Tail V@R
is weakly compact provided that $\la\in(0,1]$.
The determining set $\DDD_\mu$ of Weighted V@R is weakly compact
provided that $\mu$ is concentrated on $(0,1]$;
this follows from the explicit representation
of this set provided in~\cite{CD03}
(the proof can also be found in~\cite[Th.~4.73]{FS04}
or~\cite[Th.~1.53]{S05}); this can also be seen from the
representation of~$\DDD_\mu$ provided in~\cite{C05e}.

The following example shows that the condition $X\in L_s^1(\DDD)$
in Proposition~\ref{EM1} cannot be replaced by the condition
$X\in L_w^1(\DDD)$.

\begin{Example}\rm
\label{EM3}
Let $\Omega=[0,1]$ be endowed with the Lebesgue measure.
Consider $Z_n=\sqrt{n}I_{[0,1/n]}+1-1/\sqrt{n}$, $n\in\N$.
Then $Y_n:=Z_n-1\xra{L^1}0$, and therefore, the~set
$$
\DDD=\biggl\{1+\sum_{n=1}^\infty a_nY_n:a_n\ge0,\;
\sum_{n=1}^\infty a_n\le1\biggr\}
$$
is convex, $L^1$-closed, and uniformly integrable.
Thus, $\DDD$ is weakly compact.
Now, consider $X(\omega)=-1/\sqrt{\omega}$. Then
$\EE_\PP Z_nX=-4+2/\sqrt{n}$.
Thus, $\inf_{\QQ\in\DDD}\EE_\QQ X=-4$, while there exists
no $\QQ\in\DDD$ such that $\EE_\QQ X=-4$.\End
\end{Example}

\subsection{Capital Allocation}
\label{CA}

Let $(\Omega,\F,\PP)$ be a probability space,
$u$ be a coherent utility function with the determining
set~$\DDD$,
and let $X^1,\dots,X^d\in L_w^1(\DDD)$ be the discounted
P\&Ls produced by different components of a firm
(P\&L means the Profit\&Loss, i.e. the difference
between the terminal wealth and the initial wealth).
We will use the notation $X=(X^1,\dots,X^d)$.

Informally, the capital allocation problem is the following.
How is the total risk $\rho\bigl(\sum_i X^i\bigr)$ being split
between the components $1,\dots,d$?
In other words, we are looking for a vector $(x^1,\dots,x^d)$
such that $x^i$ means that part of the risk carried by
the $i$-th component.
Taking $x^i=\rho(X^i)$ does not work because
$\sum_i\rho(X^i)\ne\rho\bigl(\sum_i X^i\bigr)$.
The following definition of a capital allocation is
taken from~\cite[Sect.~9]{D05}.
In fact, it is closely connected with the coalitional
games (see~\cite{D01}).

\skm
\textbf{Problem (capital allocation):}
Find $x^1,\dots,x^d\in\R$ such that
\begin{gather}
\label{ca1}
\sum_{i=1}^d x^i=u\Bigl(\sum_{i=1}^d X^i\Bigr),\\
\label{ca2}
\forall h^1,\dots,h^d\in\R_+,\;\;
\sum_{i=1}^d h^i x^i\ge
u\Bigl(\sum_{i=1}^d h^i X^i\Bigr).
\end{gather}
We will call a solution of this problem a
\textit{utility allocation between $X^1,\dots,X^d$}.
A \textit{capital allocation} is defined as a utility
allocation with the minus sign.

\skm
From the financial point of view, $-x^i$ is the contribution
of the $i$-th component to the total risk of the firm, or,
equivalently, the capital that should be allocated to this
component.
In order to illustrate the meaning of~\eqref{ca2}, consider
the example $h^i=I(i\in J)$, where $J$ is a subset of
$\{1,\dots,d\}$. Then~\eqref{ca2} means that the capital
allocated to a part of the firm does not exceed the risk
carried by that part.

\skm
Let us introduce the notation
$G=\cl\{\EE_\QQ X:\QQ\in\DDD\}$, where
``$\cl$'' denotes the closure.
Note that $G$ is convex and compact.
We will call it the \textit{generating set}
or simply the \textit{generator} for~$X$ and~$u$.
This term is justified by the line
\begin{equation}
\label{ca2.5}
u(\lb h,X\rb)
=\inf_{\QQ\in\DDD}\EE_\QQ\lb h,X\rb
=\inf_{\QQ\in\DDD}\lb h,\EE_\QQ X\rb
=\min_{x\in G}\lb h,x\rb,\quad h\in\R^d.
\end{equation}
Note that the last expression is a classical object of
convex analysis known as the \textit{support function}
of the convex set~$G$.

\begin{Theorem}
\label{CA1}
The set $U$ of utility allocations between
$X^1,\dots,X^d$ has the form
\begin{equation}
\label{ca3}
U=\argmin_{x\in G}\lb e,x\rb,
\end{equation}
where $e=(1,\dots,1)$.
Furthermore, for any utility allocation $x$, we have
\begin{equation}
\label{ca4}
\forall h^1,\dots,h^d\in\R,\;\;
\sum_{i=1}^d h^i x^i\ge
u\Bigl(\sum_{i=1}^d h^i X^i\Bigr)
\end{equation}
If moreover $X^1,\dots,X^d\in L_s^1(\DDD)$ and $\DDD$ is weakly
compact, then
\begin{equation}
\label{ca5}
U=\biggl\{\EE_\QQ X:\QQ\in\EX_\DDD\Bigl(\sum_{i=1}^dX^i\Bigr)\biggr\}.
\end{equation}
\end{Theorem}

\sksmin
\Proof
(The proof is illustrated by Figure~1.)
For $h\in\R^d$, we set
\begin{align*}
L(h)&=\Bigl\{x\in\R^d:\lb h,x\rb=\min_{y\in G}\lb h,y\rb\Bigr\},\\
M(h)&=\Bigl\{x\in\R^d:\lb h,x\rb\ge\min_{y\in G}\lb h,y\rb\Bigr\}.
\end{align*}
It is seen from~\eqref{ca2.5} that the set of points
$x\in\R^d$ that satisfy~\eqref{ca1} is $L(e)$.
The set of points~$x$ that satisfy~\eqref{ca2} is
$\bigcap_{h\in\R_+^d}M(h)=G+\R_+^d$.
The set of points~$x$ that satisfy~\eqref{ca4} is
$\bigcap_{h\in\R^d}M(h)=G$.
This proves~\eqref{ca3} and~\eqref{ca4}.
Furthermore, the set $\{\EE_\QQ X:\QQ\in\DDD\}$ is closed
(the proof is similar to the proof of Proposition~\ref{EM2}).
Now, equality~\eqref{ca5} follows immediately from~\eqref{ca3}
and the definition of $\EX_\DDD$.\End

\begin{figure}[h]
\begin{picture}(150,70)(-41,-12.5)
\put(-0.3,-3.5){\includegraphics{poem.1}}
\put(16.9,20){\small $M(h)$}
\put(36.2,17){\small $U$}
\put(5.9,15){\small $L(h)$}
\put(3,39){\small $L(e)$}
\put(60,41){\small $G$}
\put(60.3,50){\small $G+\R_+^d$}
\put(7,5){\small $e$}
\put(-2,9){\small $h$}
\put(-6,-10){\small\textbf{Figure~1.} Solution of the capital
allocation problem}
\end{picture}
\end{figure}

\skm
If $G$ is strictly convex (i.e. its interior is
nonempty and its border contains no interval),
then a utility allocation is unique.
However, in general it is not unique as shown by the
example below.

\begin{Example}\rm
\label{CA2}
Let $d=2$ and $X^2=-X^1$. Then $G$ is the interval with
the endpoints $(u(X^1),-u(X^1))$ and $(-u(-X^1),u(-X^1))$.
In this example, $U=G$.\End
\end{Example}

Let us now find the solution of the capital allocation
problem in the Gaussian case.

\begin{Example}\rm
\label{CA3}
Let $X$ have Gaussian distribution with mean $a$
and covariance matrix~$C$.
Let~$u$ be a \textit{law invariant} coherent utility
function, i.e. $u(X)$ depends only on the distribution of~$X$;
we also assume that~$u$ is finite on Gaussian
random variables.

Then there exists $\ga>0$ such that, for a Gaussian random
variable~$\xi$ with mean~$m$ and variance~$\si^2$, we have
$u(\xi)=m-\ga\si$.
Let $L$ denote the image of $\R^d$ under the map
$x\mapsto Cx$. Then the inverse $C^{-1}:L\to L$
is correctly defined.
It is easy to see that
$$
G=a+\{C^{1/2}x:\|x\|\le\ga\}
=a+\{y\in L:\lb y,C^{-1}y\rb\le\ga^2\}.
$$

Let $e=(1,\dots,1)$ and assume first that $Ce\ne0$.
In this case the utility allocation~$x_0$ between
$X^1,\dots,X^d$ is determined uniquely.
In order to find it, note that, for any $y\in L$ such that
$$
\frac{d}{d\eps}\Bigl|_{\eps=0}\lb x_0-a+\eps y,C^{-1}(x_0-a+\eps y)\rb=0,
$$
we have $\lb e,y\rb=0$.
This implies that $C^{-1}(x_0-a)=\al\pr_L e$ with some
constant~$\al$ ($\pr_L$ denotes the orthogonal projection
on~$L$).
Thus, $x_0=a+\al Ce$.
As $x_0$ should belong to the relative border of~$G$
(i.e. the border in the relative topology of $a+L$), we
have $\lb x_0-a,C^{-1}(x_0-a)\rb=\ga^2$, i.e.
$\al=-\ga\lb e,Ce\rb^{-1/2}$.
As a result, the utility allocation between $X^1,\dots,X^d$
is $a-\ga\lb e,Ce\rb^{-1/2}Ce$.

Assume now that $Ce=0$. This means that $e$ is orthogonal
to~$L$, and then the set of utility allocations between
$X^1,\dots,X^d$ is~$G$.

Let us remark that in this example the solution of the
capital allocation problem depends on~$u$ rather weakly,
i.e. it depends only on~$\ga$.\End
\end{Example}

\subsection{Risk Contribution}
\label{RC}

Let $(\Omega,\F,\PP)$ be a probability space,
$u$ be a coherent utility function with the determining set~$\DDD$,
$X\in L^0$ be the discounted P\&L produced by
a component of some firm, and $Y\in L^0$ be the
discounted P\&L produced by the whole firm.

From the financial point of view, such a firm assesses
the risk of~$X$ not as $\rho(X)$ but rather as
$\rho(W+X)-\rho(X)$.
Below we define a risk contribution $\rho^c(X;W)$
in such a way that it is a coherent risk measure as a
function of~$X$ and $\rho^c(X;W)\approx\rho(W+X)-\rho(W)$
provided that $X$ is small as compared to~$W$
(the precise statement is Theorem~\ref{RC2}).

\begin{Definition}\rm
\label{RC1}
The \textit{utility contribution of $X$ to $Y$} is
$$
u^c(X;W)=\inf_{\QQ\in\EX_\DDD(Y)}\EE_\QQ X.
$$
The \textit{risk contribution of $X$ to $Y$} is defined as
$\rho^c(X;Y)=-u^c(X;Y)$.
\end{Definition}

The utility contribution is a coherent utility function
provided that $\EX_\DDD(Y)\ne\emp$.

If $\DDD$ is weakly compact and $X,Y\in L_s^1(\DDD)$
then, by Theorem~\ref{CA1},
$$
u^c(X;Y)=\inf\{x^1:(x^1,x^2)\text{ is a utility allocation
between }X,Y-X\}.
$$
This formula enables one to define risk contribution under
a weaker assumption $X,Y\in L_w^1(\DDD)$.

If $\DDD$ is weakly compact, $X^1,\dots,X^d\in L_s^1(\DDD)$,
and $\EX_\DDD\bigl(\sum_i X^i\bigr)$ is a singleton,
then (in view of Theorem~\ref{CA1}) the utility
allocation between $X^1,\dots,X^d$ is unique and has the form
$$
\biggl(u^c\Bigl(X^1;\sum_{i=1}^d X^i\Bigr),\dots,
u^c\Bigl(X^d;\sum_{i=1}^d X^i\Bigr)\biggr).
$$
This shows the relevance of the given definition.
Another argument supporting this definition is the statement below.

\begin{Theorem}
\label{RC2}
If $\DDD$ is weakly compact and $X,Y\in L_s^1(\DDD)$, then
$$
u^c(X;Y)=\lim_{\eps\da0}\eps^{-1}(u(Y+\eps X)-u(Y)).
$$
\end{Theorem}

\sksmin
\Proof
(The proof is illustrated by Figure~2.)
Consider the generator $G=\cl\{\EE_\QQ(X,Y):\QQ\in\DDD\}$
and set
$b=\inf\{y:\exists x:(x,y)\in G\}$,
$I=\{x:(x,b)\in G\}$,
$J=\{x:\exists y:(x,y)\in G\}$,
$a=\inf\{x:x\in I\}$.
Note that $u^c(X;Y)=a$.
The minimum $\min_{(x,y)\in G}\lb(\eps,1),(x,y)\rb$ is
attained at a point $(a(\eps),b(\eps))$.
We have $a(\eps)\le a$, $b(\eps)\ge b$, and
$(a(\eps),b(\eps))\xra[\eps\da0]{}(a,b)$.
Furthermore,
$\eps a(\eps)+b(\eps)\le\eps a+b$,
which implies that
$0\le b(\eps)-b\le\eps(a-a(\eps))$.
As a result,
\begin{align*}
\lim_{\eps\da0}\eps^{-1}(u(Y+\eps X)-u(Y))
&=\lim_{\eps\da0}\eps^{-1}(\eps a(\eps)+b(\eps)-b)\\
&=a+\lim_{\eps\da0}\eps^{-1}(b(\eps)-b)
=a
=u^c(X;Y).
\end{align*}

\begin{figure}
\begin{picture}(150,55)(-55,-15)
\put(4.7,-0.5){\includegraphics{poem.2}}
\put(0,0){\vector(1,0){50}}
\put(0,0){\vector(0,1){40}}
\multiput(0,5)(2,0){10}{\line(1,0){1}}
\multiput(5,0)(0,2){10}{\line(0,1){1}}
\multiput(20,0)(0,2){3}{\line(0,1){1}}
\multiput(30,0)(0,2){3}{\line(0,1){1}}
\multiput(45,0)(0,2){10}{\line(0,1){1}}
\put(-3.5,38){\small $y$}
\put(-3.5,4){\small $b$}
\put(19,-4){\small $a$}
\put(24,-4){\small $I$}
\put(36,-4){\small $J$}
\put(48,-4){\small $x$}
\put(23,19){\small $G$}
\put(17,-15){\small\textbf{Figure~2}}
\end{picture}
\end{figure}

\begin{Example}\rm
\label{RC3}
{\bf(i)} Let $Y$ be a constant.
In this case $\EX_\DDD(Y)=\DDD$, so that
$u^c(X;Y)=u(X)$.

{\bf(ii)} Let $X=\al Y$ with $\al\in\R_+$.
Then $u^c(X;Y)=\al u(Y)$.

{\bf(iii)} Let $X,Y$ have a jointly Gaussian distribution
with mean $(\EE X,\EE Y)$ and covariance matrix~$C$.
Let $u$ be a law invariant coherent utility function
that is finite on Gaussian random variables.
Then there exists $\ga>0$ such that, for a Gaussian random
variable~$\xi$ with mean~$m$ and variance~$\si^2$, we have
$u(\xi)=m-\ga\si$.
Assume that $X$ and $Y$ are not degenerate and
$\corr(X,Y)\ne\pm1$. It follows from Example~\ref{CA3} that
\begin{align*}
u^c(X;Y)
&=\EE X-\gamma\lb e_2,C e_2\rb^{-1/2}Ce_2\\
&=\EE X-\ga\,\frac{\cov(X,Y)}{(\var Y)^{1/2}}\\
&=\EE X+(u(X)-\EE X)\corr(X,Y),
\end{align*}
where $e_2=(0,1)$. In particular, if $\EE X=\EE Y=0$, then
$$
\frac{u^c(X;Y)}{u(X)}
=\corr(X,Y)
=\frac{\text{V@R}^c(X;Y)}{\text{V@R}(X)},
$$
where $\var$ denotes the variance and
$\text{V@R}^c$ denotes the V@R contribution
(for the definition, see~\cite[Sect.~7]{M02}).\End
\end{Example}

\section{Good Deals Pricing}
\label{P}

\subsection{Utility-Based Good Deals Pricing}
\label{UP}

Let $(\Omega,\F,\PP)$ be a probability space,
$u$ be a coherent utility function with the weakly compact
determining set~$\DDD$,
and $A$ be a convex subset of~$L^0$.
From the financial point of view, $A$ is the set of
various discounted P\&Ls that can be obtained in the model
under consideration by employing various trading strategies
(examples are given in Subsections~\ref{PF}--\ref{PDT}).
It will be called the \textit{set of attainable P\&Ls}.
We will assume that $A$ is $\DDD$-consistent
(see Definition~\ref{UP2} below).
It is shown in Subsections~\ref{PF}--\ref{PDT} that
this assumption is automatically satisfied for natural models.

First, we give the definition of a risk-neutral measure.
Of course, this notion is a classical object of financial
mathematics, but the particular definition we need is
taken from~\cite{C05a} (it is adapted to the $L^0$-case).

\begin{Definition}\rm
\label{UP1}
A \textit{risk-neutral measure} is a measure
$\QQ\in\PPP$ such that $\EE_\QQ X\le0$ for any $X\in A$
(we use the convention $\EE X=\EE X^+-\EE X^-$,
$\infty-\infty=-\infty$).

The set of risk-neutral measures will be
denoted by $\RRR$ or by $\RRR(A)$ if there is a risk of
ambiguity.
\end{Definition}

\begin{Definition}\rm
\label{UP2}
We will say that \textit{$A$ is $\DDD$-consistent} if
there exists a set $A'\subseteq A\cap L_s^1(\DDD)$ such that
$\DDD\cap\RRR=\DDD\cap\RRR(A')$.
\end{Definition}

\begin{Definition}\rm
\label{UP3}
A model satisfies the \textit{utility-based NGD} condition
if there exists no $X\in A$ such that $u(X)>0$.
\end{Definition}

\begin{Theorem}[Fundamental Theorem of Asset Pricing]
\label{UP4}
A model satisfies the NGD condition if and only if
$\DDD\cap\RRR\ne\emp$.
\end{Theorem}

\sksminus
\Proof
The ``if'' part is obvious.
Let us prove the ``only if'' part.

Fix $X_1,\dots,X_M\in A'$.
It follows from the weak continuity of the
maps $\DDD\ni\QQ\mapsto\EE_\QQ X_m$ that the set
$G=\{\EE_\QQ(X_1,\dots,X_M):\QQ\in\DDD\}$ is compact.
Clearly, $G$ is convex.
Suppose that $G\cap(-\infty,0]^M=\emp$.
Then there exist $h\in\R^M$ and $\eps>0$ such that
$\lb h,x\rb\ge\eps$ for any $x\in G$
and $\lb h,x\rb\le0$ for any $x\in(-\infty,0]^M$.
Hence, $h\in\R_+^M$.
Without loss of generality, $\sum_m h_m=1$.
Then $X=\sum_m h_m X_m\in A$ and
$\EE_\QQ X\ge\eps$ for any $\QQ\in\DDD$,
so that $u(X)>0$.

The obtained contradiction shows that,
for any $X_1,\dots,X_M\in A'$, the set
$$
B(X_1,\dots,X_M)=\{\QQ\in\DDD:\EE_\QQ X_m\le0
\text{ for any }m=1,\dots,M\}
$$
is nonempty. As $X_m\in L_s^1(\DDD)$, the map
$\DDD\ni\QQ\mapsto\EE_\QQ X_m$ is
weakly continuous, and therefore, $B(X_1,\dots,X_M)$ is
weakly closed.
Furthermore, any finite intersection of sets of this form
is nonempty.
Consequently, there exists a measure $\QQ$ that belongs
to each~$B$.
Then $\EE_\QQ X\le0$ for any $X\in A'$,
which means that $\QQ\in\DDD\cap\RRR(A')$.
As $A$ is $\DDD$-consistent, $\QQ\in\DDD\cap\RRR$.\End

\skm
\Remarks
\texttt{(i)} As opposed to the fundamental theorems of asset
pricing dealing with the NA condition
and its strengthenings (see~\cite{C05a}, \cite{DS94},
\cite{DS98}), here we need not take any closure
of~$A$ when defining the NGD.
Essentially, this is the compactness of~$\DDD$ that yields
the fundamental theorem of asset pricing.

\texttt{(ii)} If $\DDD=\PPP$, then the NGD condition
means that there exists no $X\in A$ with
$\essinf_\omega X(\omega)>0$.
This is very close to the NA condition.
However, in this case $\DDD$ is not uniformly integrable
and Theorem~\ref{UP4} might be violated.
Indeed, let $A=\{hX:h\in\R\}$, where $X$
has uniform distribution on $[0,1]$.
Then the NGD is satisfied, while $\RRR=\emp$.

\skm
Now, let $F\in L^0$ be the discounted payoff of a contingent claim.

\begin{Definition}\rm
\label{UP5}
A \textit{utility-based NGD price} of~$F$ is a real number~$x$ such that
the extended model $(\Omega,\F,\PP,\DDD,A+\{h(F-x):h\in\R\})$
satisfies the NGD condition.

The set of the NGD prices will be denoted by $\INGD(F)$.
\end{Definition}

\begin{Corollary}[Fair price interval]
\label{UP6}
For $F\in L_s^1(\DDD)$,
$$
\INGD(F)=\{\EE_\QQ F:\QQ\in\DDD\cap\RRR\}.
$$
\end{Corollary}

\sksmin
\Proof
Denote $\{h(F-x):h\in\R\}$ by~$A(x)$.
Clearly, $A+A(x)$ is $\DDD$-consistent
(in order to prove this, it is sufficient to consider
$A'+A(x)$).
It follows from Theorem~\ref{UP4} that $x\in\INGD(F)$
if and only if $\DDD\cap\RRR(A+A(x))\ne\emp$.
It is easy to check that $\QQ\in\RRR(A+A(x))$
if and only if $\QQ\in\RRR$ and $\EE_\QQ F=x$.
This completes the proof.\End

\skm
\Remark
As opposed to the NA price intervals, the NGD price
intervals are closed (this follows from the
weak continuity of the map $\DDD\cap\RRR\mapsto\EE_\QQ F$).

\skm
To conclude the subsection, we will discuss the origin of~$\DDD$.
First of all, $\DDD$ might be the determining set of
a coherent utility function like Tail V@R or Weighted V@R.
The set~$\DDD$ might also correspond to a weighted average
or the minimum of several coherent utility functions.
It is also possible that $\DDD$ originates from the
classical utility maximization as described by the example
below.

\begin{Example}\rm
\label{UP7}
Let $\PP_1,\dots,\PP_N$ be a family of probability
measures, $u_1,\dots,u_N$ be a family of classical
utility functions (i.e. smooth concave increasing functions
$\R\to\R$), and $W_1,\dots,W_N$ be a family of
random variables.
From the financial point of view, $\PP_n$, $u_n$, and
$W_n$ are the subjective probability, the utility function,
and the future wealth of the $n$-th market participant,
respectively. Consider a measure $\QQ_n=c_n u'_n(W_n)\PP_n$,
where $c_n$ is the normalizing constant.
Then, for any trading opportunity $X\in L^0$, we have
\begin{equation}
\label{up1}
\frac{d}{d\eps}\Bigl|_{\eps=0}u_n(W_n+\eps X)
=\EE_{\PP_n}u'_n(W_n)X
=\EE_{\QQ_n}c_n^{-1} X
\end{equation}
(we assume that all the expectations exist and
integration is interchangeable with differentiation).
Thus, an opportunity $\eps X$ with a small $\eps>0$
is attractive to the
$n$-th participant if and only if $\EE_{\QQ_n}X>0$, so that
$\QQ_n$ might be called the valuation measure
of the $n$-th participant.
Take $\DDD=\conv(\QQ_1,\dots,\QQ_N)$
and consider the corresponding coherent utility function~$u$.
Then $u(X)>0$ if and only if $\EE_{\QQ_n}X>0$ for any $n$.
In view of~\eqref{up1}, this means that $\eps X$
with some $\eps>0$ is attractive to any market participant
(this is similar to the notion of a strictly
acceptable opportunity introduced in~\cite{CGM01}).
Thus, in this example the NGD means the absence of a trading
opportunity that is attractive to every agent.\End
\end{Example}

\subsection{RAROC-Based Good Deals Pricing}
\label{RP}

Let $(\Omega,\F,\PP)$ be a probability space,
$\RD\subset\PPP$ be a convex weakly compact set,
$\PD$ be an $L^1$-closed convex subset of~$\RD$,
and $A$ be a convex subset of~$L^0$.
We will call $\PD$ the \textit{profit-determining set}.
Thus, the profit of a position that yields a P\&L~$X$
is $\inf_{\QQ\in\PD}\EE_\QQ X$.
We will call $\RD$ the \textit{risk-determining set},
so that the risk of a position that yields a P\&L~$X$
is $-\inf_{\QQ\in\DDD}\EE_\QQ X$.
A canonical example is:
$\PD=\{\PP\}$ and $\RD$ is the determining set of
a coherent utility function.
We will assume that $A$ is $\RD$-consistent.
Finally, we fix a positive number~$R$
meaning the upper limit on a possible RAROC.

\begin{Definition}\rm
\label{RP1}
The \textit{Risk-Adjusted Return on Capital} (RAROC) for
$X\in L^0$ is defined as
$$
\RAROC(X)=\begin{cases}
+\infty&\text{if}\;\;\inf_{\QQ\in\PD}\EE_\QQ X>0
\text{ and }\inf_{\QQ\in\RD}\EE_\QQ X\ge0,\\[2mm]
\ds\frac{\inf_{\QQ\in\PD}\EE_\QQ X}{-\inf_{\QQ\in\RD}\EE_\QQ X}
&\text{otherwise}
\end{cases}
$$
with the convention $\frac{0}{0}=0$, $\frac{\infty}{\infty}=0$.
\end{Definition}

\begin{Definition}\rm
\label{RP2}
A model satisfies the \textit{RAROC-based NGD} condition
if there exists no $X\in A$ such that
$\RAROC(X)>R$.
\end{Definition}

\begin{Theorem}[Fundamental Theorem of Asset Pricing]
\label{RP3}
A model satisfies the NGD condition if and only if
\begin{equation}
\label{rp1}
\biggl(\frac{1}{1+R}\,\PD+\frac{R}{1+R}\,\RD\biggr)
\cap\RRR\ne\emp.
\end{equation}
\end{Theorem}

\sksminus
\Proof
Let us first consider the case $R>0$.
Then, for any $X\in L^0$,
$$
\RAROC(X)>R
\;\Longleftrightarrow\;
\inf_{\QQ\in\PD}\EE_\QQ X+R\inf_{\QQ\in\RD}\EE_\QQ X>0
\;\Longleftrightarrow\;
\inf_{\QQ\in\DDD}\EE_\QQ X>0,
$$
where $\DDD=\bigl(\frac{1}{1+R}\,\PD+\frac{R}{1+R}\,\RD\bigr)$.
Clearly, $\DDD$ is weakly compact
(note that $\DDD\subseteq\RD$,
while $L_s^1(\DDD)=L_s^1(\RD)$)
and $A$ is $\DDD$-consistent.
Now, the statement follows from Theorem~\ref{UP4}.

Let us now consider the case $R=0$.
Then the ``if'' part is obvious, and we should check
the ``only if'' part.
Take $A'\subseteq A\cap L_s^1(\RD)$ such that
$\RD\cap\RRR=\RD\cap\RRR(A')$.
For any $X\in\conv A'$, $\inf_{\QQ\in\PD}\EE_\QQ X\le0$.
Repeating the arguments from the proof of Theorem~\ref{UP4},
we get $\PD\cap\RRR\ne\emp$.\End

\begin{Definition}\rm
\label{RP4}
A \textit{RAROC-based NGD price} of a contingent
claim~$F$ is a real number~$x$ such that
the extended model $(\Omega,\F,\PP,\PD,\RD,A+\{h(F-x):h\in\R\})$
satisfies the NGD condition.

The set of the NGD prices will be denoted by $\INGD(F)$.
\end{Definition}

\begin{Corollary}[Fair price interval]
\label{RP5}
For $F\in L_s^1(\DDD)$,
$$
\INGD(F)=\biggl\{\EE_\QQ F:\QQ\in\Bigl(\frac{1}{1+R}\,\PD
+\frac{R}{1+R}\,\RD\Bigr)\cap\RRR\biggr\}.
$$
\end{Corollary}

\sksmin
This statement follows from Theorem~\ref{RP3}.

\subsection{Static Model with a Finite Number of Assets}
\label{PF}

We consider the model of the previous subsection with
$A=\{\lb h,S_1-S_0\rb:h\in\R^d\}$, where $S_0\in\R^d$
and $S_1^1,\dots,S_1^d\in L_s^1(\RD)$.
From the financial point of view, $S_n^i$ is the discounted
price of the $i$-th asset at time~$n$.
Clearly, in this model $A$ is $\RD$-consistent
and $\RD\cap\RRR=\RD\cap\MMM$, where $\MMM$ is the set of
martingale measures:
$$
\MMM=\{\QQ\in\PPP:\EE_\QQ|S_1|<\infty\text{ and }
\EE_\QQ S_1=S_0\}.
$$

\skm
\Remark
We have $\MMM\subseteq\RRR$, but the reverse inclusion
might be violated.
Indeed, let $d=1$ and let $S_1$ be such that
$\EE_\PP S_1^+=\EE_\PP S_1^-=\infty$.
Then $\PP\in\RRR$, while $\PP\notin\MMM$.

\skm
Let us now provide a geometric interpretation of
Theorems~\ref{UP4} and~\ref{RP3}.
For this, we only assume that $\PD\subseteq\RD\subseteq\PPP$
are convex sets and $S_1\in L_w^1(\RD)$.
Let us introduce the notation (see Figure~3)
\begin{align*}
E&=\cl\{\EE_\QQ S_1:\QQ\in\PD\},\\
G&=\cl\{\EE_\QQ S_1:\QQ\in\RD\},\\
G_R&=\frac{1}{1+R}\,E+\frac{R}{1+R}\,G,\\
D&=\conv\supp\Law_\PP S_1,
\end{align*}
where ``$\supp$'' denotes the support,
and let $D^\circ$ denote the relative interior of~$D$
(i.e. the interior in the relative topology of the smallest
affine subspace containing~$D$).
It is easy to see from the equalities
\begin{align*}
\inf_{\QQ\in\PD}\EE_\QQ\lb h,S_1-S_0\rb
&=\inf_{x\in E}\lb h,x-S_0\rb,\\
\inf_{\QQ\in\RD}\EE_\QQ\lb h,S_1-S_0\rb
&=\inf_{x\in G}\lb h,x-S_0\rb
\end{align*}
that the following equivalences are true:
\begin{align*}
\text{RAROC-based NGD}\;&\Longleftrightarrow\;S_0\in G_R,\\
\text{utility-based NGD corresponding to~$u$}
\;&\Longleftrightarrow\;S_0\in G,\\
\text{NA}\;&\Longleftrightarrow\;S_0\in D^\circ
\end{align*}
(the last equivalence is a well-known result of arbitrage
pricing; see~\cite[Ch.~V, \S~2e]{S99}).

Now, let $F\in L_w^1(\RD)$ be the discounted payoff of a
contingent claim.
Let $\wt E$, $\wt G$, $\wt G_R$,
$\wt D$, and $\wt D^\circ$ denote the versions of the sets
$E$, $G$, $G_R$, $D$, and $D^\circ$ defined for
$\wt S_1=(S_1^1,\dots,S_1^d,F)$ instead of~$S_1$.
Let $\INGDR(F)$ denote the RAROC-based NGD price interval,
$\INGD(F)$ denote the utility-based NGD price interval
(corresponding to~$u$), and $\INA(F)$ denote the NA price
interval. Then
\begin{align*}
\INGDR(F)&=\{x:(S_0,x)\in\wt G_R\},\\
\INGD(F)&=\{x:(S_0,x)\in\wt G\},\\
\INA(F)&=\{x:(S_0,x)\in\wt D^\circ\}.
\end{align*}

\begin{figure}[h]
\begin{picture}(150,87)(-61.5,-27.5)
\put(-26,0){\includegraphics{poem.3}}
\put(-5,-5){\vector(1,0){65}}
\put(-5,-5){\vector(0,1){62}}
\put(20,-5){\line(0,1){60}}
\multiput(-5.5,51.5)(1.9,0){14}{\line(1,0){1}}
\multiput(-5.5,44)(1.9,0){14}{\line(1,0){1}}
\multiput(-5.5,36.5)(1.9,0){14}{\line(1,0){1}}
\multiput(-5.5,15.5)(1.9,0){14}{\line(1,0){1}}
\multiput(-5.5,8)(1.9,0){14}{\line(1,0){1}}
\multiput(-5.5,0.5)(1.9,0){14}{\line(1,0){1}}
\put(57,-9){\small $\R^d$}
\put(-9,55){\small $\R$}
\put(18,-9){\small $S_0$}
\put(24,26){\small $\wt E$}
\put(23,33){\small $\wt G_R$}
\put(24,40){\small $\wt G$}
\put(23.3,47){\small $\wt D^\circ$}
\put(-11,34){\scalebox{0.8}{\rotatebox{270}{\small $\INGDR(F)$}}}
\put(-19,32){\scalebox{0.8}{\rotatebox{270}{\small $\INGD(F)$}}}
\put(-27,31.5){\scalebox{0.8}{\rotatebox{270}{\small $\INA(F)$}}}
\put(-19,-22){\parbox{75mm}{\small\textbf{Figure~3.}
The geometric representation of price intervals provided
by various techniques}}
\end{picture}
\end{figure}

\begin{Example}\rm
\label{PF1}
Let $S_1$ have Gaussian distribution with mean $a$ and
covariance matrix~$C$.
Let $\PD=\{\PP\}$ and $\RD$ be the determining set of a
law invariant coherent utility function~$u$ that
is finite on Gaussian random variables.
Let $F$ be such that the vector
$(S_1^1,\dots,S_1^d,F)$ is Gaussian.
Denote $c=\cov(S_1,F)$ (we use the vector form
of notation).

There exists $b\in\R^d$ such that
$Cb=c$. We can write $F=\lb b,S_1-a\rb+\EE F+\wt F$.
Then $\EE\wt F=0$ and $\cov(\wt F,S_1)=0$, so that $\wt F$ is
independent of $S_1$.
Note that
$$
\si^2:=\var\wt F
=\var F-\var\lb b,S_1-a\rb
=\var F-\lb b,Cb\rb
=\var F-\lb b,c\rb.
$$
Clearly, if $\si^2=0$, then
$$
\INGDR(F)=\INGD(F)=\INA(F)=\{\lb b,S_0-a\rb+\EE F\}.
$$
Let us now assume that $\si^2>0$.

Obviously, $\INA(F)=\R$.

In order to find $\INGD(F)$, note that
$\INGD(F)=\lb b,S_0-a\rb+\EE F+\INGD(\wt F)$.
Let $L$ denote the image of $\R^d$ under the map
$x\mapsto Cx$. Then the inverse $C^{-1}:L\to L$ is correctly
defined.
As $u$ is law invariant, there exists $\ga>0$ such that,
for a Gaussian random variable $\xi$ with mean~$m$ and
variance~$\si^2$, we have $u(\xi)=m-\ga\si$.
From this, it is easy to see that the set
$\wt G:=\{\EE_\QQ(S_1,\wt F):\QQ\in\RD\}$ has the form
\begin{equation}
\label{pf0}
\wt G=(a,0)+\{(x,y):x\in L,\;y\in\R:\lb x,C^{-1}x\rb+\si^{-2}y^2\le\ga^2\}.
\end{equation}
Consequently,
$$
\INGD(F)
=\bigl[\lb b,S_0-a\rb+\EE F-\al,\lb b,S_0-a\rb+\EE F+\al\bigr],
$$
where $\al=(\si^2\ga^2-\si^2\lb S_0-a,C^{-1}(S_0-a)\rb)^{1/2}$.
(In particular, the NGD is satisfied if and only if
$\lb S_0-a,C^{-1}(S_0-a)\rb\le\ga^2$.)

Similar arguments show that
$$
\INGDR(F)
=\bigl[\lb b,S_0-a\rb+\EE F-\al(R),\lb b,S_0-a\rb+\EE F+\al(R)\bigr],
$$
where $\al(R)=\bigl(\frac{\si^2\ga^2R^2}{(1+R)^2}-\si^2\lb S_0-a,C^{-1}(S_0-a)\rb\bigr)^{1/2}$.
(In particular, the NGD({\sl R}) condition is satisfied
if and only if
$\lb S_0-a,C^{-1}(S_0-a)\rb\le\frac{\ga^2R^2}{(1+R)^2}$.)

Let us remark that $\INGD(F)$ and $\INGDR(F)$ depend on $u$
rather weakly, i.e. they depend only on~$\ga$.\End
\end{Example}

\subsection{Dynamic Model with an Infinite Number of Assets}
\label{PDF}

Let $(\Omega,\F,(\F_t)_{t\in[0,T]},\PP)$ be a filtered
probability space. We assume that $\F_0$ is trivial.
Let $\DDD\subseteq\PPP$ be a convex weakly compact set
(in the framework of Subsection~\ref{UP}, $\DDD$ is the
determining set of~$u$; in the framework of Subsection~\ref{RP},
$\DDD=\frac{1}{1+R}\,\PD+\frac{R}{1+R}\,\RD$).
Let $(S^i),\;i\in I$ be a family of $(\F_t)$-adapted
c\`adl\`ag processes (the set $I$ is arbitrary and we
impose no assumptions on the probabilistic structure of~$S^i$
like the assumption that $S^i$ is a semimartingale).
From the financial point of view, $S^i$ is the discounted
price process of the $i$-th asset.
We assume that $S_t^i\in L_s^1(\DDD)$ for any $t\in[0,T]$,
$i\in I$.
The set of P\&Ls an agent can obtain by piecewise constant
trading strategies (and only such strategies can be employed
in practice) is naturally defined as
\begin{equation}
\begin{split}
\label{pdf1}
A&=\biggl\{\sum_{n=1}^N\sum_{i\in I}
H_n^i(S_{u_n}^i-S_{u_{n-1}}^i):
N\in\N,\;u_0\le\dots\le u_N,
\text{ are }(\F_t)\text{-stopping times, }\\
&\hspace*{7mm}H_n^i\text{ is }\F_{u_{n-1}}\text{-measurable, and }
H_n^i=0\text{ for all }i,\text{ except for a finite set}\biggr\}.
\end{split}
\end{equation}

\begin{Lemma}
\label{PDF1}
We have
$\DDD\cap\RRR=\DDD\cap\RRR(A')=\DDD\cap\MMM$,
where
\begin{align*}
A'&=\{H(S_v^i-S_u^i):u\le v\in[0,T],\;i\in I,\;H\;
\text{\rm is }\F_u\text{\rm-measurable and bounded}\},\\
\MMM&=\{\QQ\in\PPP:\text{\rm for any }i\in I,\;S^i\;
\text{\rm is an }(\F_t,\QQ)\text{\rm-martingale}\}.
\end{align*}
\end{Lemma}

\sksminus
\Proof
The inclusions
$\DDD\cap\RRR\subseteq\DDD\cap\RRR(A')\subseteq\DDD\cap\MMM$
are clear. So, it is sufficient to prove the inclusion
$\DDD\cap\MMM\subseteq\DDD\cap\RRR$.
Let $\QQ\in\DDD\cap\MMM$. Take
$X=\sum_{n=1}^N\sum_{i\in I}H_n^i(S_{u_n}^i-S_{u_{n-1}}^i)\in A$.
The process
$$
M_k=\sum_{n=1}^k\sum_{i\in I}
H_n^i(S_{u_n}^i-S_{u_{n-1}}^i),\quad k=0,\dots,N
$$
is an $(\F_{u_k},\QQ)$-local martingale.
Suppose that $\EE_\QQ X^-<\infty$
(otherwise, $\EE_\QQ X=-\infty$).
Then $M$ is a martingale (see~\cite[Ch.~II,~\S~1c]{S99}),
and hence, $\EE_\QQ X=\EE_\QQ M_N=0$.
Thus, in any case, $\EE_\QQ X\le0$,
which proves that $\QQ\in\RRR$.\End

\begin{Example}\rm
\label{PDF2}
Let us consider the Black-Scholes model in the
framework of the RAROC-based pricing.
Thus, $S_t=S_0e^{\mu t+\si B_t}$, where $B$ is a
Brownian motion; we are given a risk-determining
set~$\RD$, and we take $\PD=\{\PP\}$.
Surprisingly enough, in this model
$\sup_{X\in A}\RAROC(X)=\infty$.
Indeed, the set $\MMM$ consists of a unique measure~$\QQ_0$
and $\frac{d\QQ_0}{d\PP}$ is not bounded away from zero,
so that condition~\eqref{rp1} is violated for any $R>0$.

Let us construct explicitly a sequence $X_n\in A$ with
$\RAROC(X_n)\to\infty$.
Consider
$D_n=\bigl\{\frac{d\QQ_0}{d\PP}<n^{-1}\bigr\}$ and set
$X_n=a_n I(D_n)-I(\Omega\setminus D_n)$, where $a_n$
is chosen in such a way that $\EE_{\QQ_0}X_n=0$.
Then $\EE_\PP X_n\to\infty$, while
$\inf_{\QQ\in\RD}\EE_\QQ X\ge-1$, so that
$\RAROC(X_n)\to\infty$.
Actually, $X_n\notin A$, but, for each $n$, there exists
a sequence $(Y_n^m)\in A$ such that $-2\le Y_n^m\le a_n+1$
and $Y_n^m\xra[m\to\infty]{\PP}X_n$ (we leave this to the
reader as an exercise).
Then $\RAROC(Y_n^m)\xra[m\to\infty]{}\RAROC(X_n)$, so that
$\RAROC\bigl(Y_n^{m(n)}\bigr)\to\infty$ for some
subsequence $m(n)$.

This example shows that complete models are typically
inconsistent with the RAROC-based NGD pricing.
But this technique is primarily aimed at incomplete
models because in complete ones the NA price intervals
are already exact.

Let us also remark that the utility-based NGD condition
might be naturally satisfied in the Black-Scholes model.\End
\end{Example}

\subsection{Dynamic Model with Transaction Costs}
\label{PDT}

Let $(\Omega,\F,(\F_t)_{t\in[0,T]},\PP)$ be a filtered
probability space. We assume that $\F_0$ is trivial and
$(\F_t)$ is right-continuous.
Let $\DDD\subseteq\PPP$ be a convex weakly compact set.
Let $S^{ai},S^{bi},\;i\in I$
be two families of $(\F_t)$-adapted c\`adl\`ag processes.
From the financial point of view, $S^{ai}$
(resp., $S^{bi})$ is the discounted ask (resp., bid)
price process of the $i$-th asset, so that
$S^a\ge S^b$ componentwise.
We assume that $S_t^{ai},S_t^{bi}\in L_s^1(\DDD)$
for any $t\in[0,T]$, $i\in I$.
The set of P\&Ls that can be obtained in this model is
naturally defined as
\begin{align*}
A&=\biggl\{\sum_{n=0}^N\sum_{i\in I}
\bigl[-H_n^iI(H_n^i>0)S_{u_n}^{ai}
-H_n^iI(H_n^i<0)S_{u_n}^{bi}\bigr]:\\
&\hspace*{7mm}N\in\N,\;u_0\le\dots\le u_N\text{ are }
(\F_t)\text{-stopping times, }
H_n^i\text{ is }\F_{u_n}\text{-measurable, }\\
&\hspace*{7mm}H_n^i=0\text{ for all }i,
\text{ except for a finite set, and }
\sum_{n=0}^N H_n^i=0\text{ for any }i\biggr\}.
\end{align*}
Here $H_n^i$ means the amount of the $i$-th asset
that is bought at time~$u_n$
(so that $\sum_{k=0}^n H_k^i$
is the total amount of the $i$-th asset held at
time~$u_n$).
Note that if there are no transaction costs, i.e.
$S^{ai}=S^{bi}=S^i$ for each~$i$,
then the set of attainable P\&Ls coincides with
the set given by~\eqref{pdf1}.

\begin{Lemma}
\label{PDT1}
We have
$\DDD\cap\RRR=\DDD\cap\RRR(A')=\DDD\cap\MMM$,
where
\begin{align*}
A'&=\{G(S_v^{bi}-S_u^{ai})+H(-S_v^{ai}+S_u^{bi}):
i\in I,\;u\le v\;\text{\rm are simple }(\F_t)\text{-}\\
&\hspace*{7mm}\text{\rm stopping times, }G,H\;\text{\rm are positive, bounded, }\F_u\text{\rm-measurable}\},\\
\MMM&=\{\QQ\in\PPP:\text{\rm for any }i,\;\text{\rm there exists an }
(\F_t,\QQ)\text{-}\\
&\hspace*{7mm}\text{\rm martingale }M^i\;\text{\rm such that }
S^{bi}\le M^i\le S^{ai}\}.
\end{align*}
{\rm(}A stopping time is simple if it takes on a finite
number of values.{\rm)}
\end{Lemma}

\sksminus
\Proof
The inclusion $\DDD\cap\RRR\subseteq\DDD\cap\RRR(A')$
is obvious.

Let us prove the inclusion $\DDD\cap\RRR(A')\subseteq\DDD\cap\MMM$.
Take $\QQ\in\DDD\cap\RRR(A')$.
Fix $i\in I$.
For any simple stopping times $u\le v$, we have
$S_u^{ai},S_u^{bi},S_v^{ai},S_v^{bi}\in L_s^1(\DDD)$ and
\begin{equation}
\label{pdt1}
\EE_\QQ(S_v^{ai}\mid\F_u)\ge S_u^{bi},\qquad
\EE_\QQ(S_v^{bi}\mid\F_u)\le S_u^{ai}.
\end{equation}
Consider the Snell envelopes
\begin{align*}
X_t&=\esssup_{\tau\in{\cal T}_t}\EE_\QQ(S_\tau^{bi}\mid\F_t),\quad
t\in[0,T],\\
Y_t&=\essinf_{\tau\in{\cal T}_t}\EE_\QQ(S_\tau^{ai}\mid\F_t),\quad
t\in[0,T],
\end{align*}
where ${\cal T}_t$ denotes the set of simple
$(\F_t)$-stopping times such that $\tau\ge t$.
(Recall that $\esssup_\al\xi_\al$ is a random variable
$\xi$ such that, for any $\al$, $\xi\ge\xi_\al$~a.s. and
for any other random variable $\xi'$ with this property,
we have $\xi\le\xi'$~a.s.)
Then $X$ is an $(\F_t)$-supermartingale, while $Y$ is an
$(\F_t,\QQ)$-submartingale (see~\cite[Th.~2.12.1]{EK81}).

Let us prove that, for any $t\in[0,T]$, $X_t\le Y_t$ $\QQ$-a.s.
Assume that there exists $t$ such that
$\PP(X_t>Y_t)>0$. Then there exist
$\tau,\sigma\in{\cal T}_t$ such that
$$
\QQ\bigl(\EE_\QQ(S_\tau^{bi}\mid\F_t)
>\EE_\QQ(S_\si^{ai}\mid\F_t)\bigr)>0.
$$
This implies that $\QQ(\xi>\eta)>0$, where
$\xi=\EE_\QQ(S_\tau^{bi}\mid\F_{\tau\wedge\si})$ and
$\eta=\EE_\QQ(S_\si^{ai}\mid\F_{\tau\wedge\si})$.
Assume first that
$\QQ(\{\xi>\eta\}\cap\{\tau\le\si\})>0$.
On the set $\{\tau\le\si\}$ we have
$$
\xi=S_\tau^{bi}=S_{\tau\wedge\si}^{bi},\qquad
\eta=\EE_\QQ(S_\si^{ai}\mid\F_{\tau\wedge\si})
=\EE_\QQ(S_{\tau\vee\si}^{ai}\mid\F_{\tau\wedge\si}),
$$
and we obtain a contradiction with~\eqref{pdt1}.
In a similar way we get a contradiction if we assume that
$\QQ(\{\xi>\eta\}\cap\{\tau\ge\si\})>0$.
As a result, $X_t\le Y_t$ $\QQ$-a.s.
Now, it follows from~\cite[Lem.~3]{JK95} that
there exists an $(\F_t,\QQ)$-martingale~$M$ such that
$X\le M\le Y$.
As a result, $\QQ\in\MMM$.

Let us prove the inclusion
$\DDD\cap\MMM\subseteq\DDD\cap\RRR$.
Take $\QQ\in\DDD\cap\MMM$, so that, for any~$i$, there
exists an $(\F_t,\QQ)$-martingale $M^i$ such that
$S^{bi}\le M^i\le S^{ai}$. For any
$$
X=\sum_{n=0}^N\sum_{i\in I}
\bigl[-H_n^iI(H_n^i>0)S_{u_n}^{ai}
-H_n^iI(H_n^i<0)S_{u_n}^{bi}\bigr]\in A,
$$
we have
$$
X\le\sum_{n=0}^N\sum_{i\in I}
\bigl[-H_n^iI(H_n^i\!>\!0)M_{u_n}^i
\!-H_n^iI(H_n^i\!<\!0)M_{u_n}^i\bigr]
=\sum_{n=1}^N\sum_{i\in I}\Bigl(\sum_{k=0}^{n-1}H_k^i\Bigr)
(M_{u_n}^i-M_{u_{n-1}}^i).
$$
Repeating the arguments used in the proof of Lemma~\ref{PDF1},
we get $\EE_\QQ X\le0$.
As a result, $\QQ\in\RRR$.\End

\skm
Consider now a model with proportional
transaction costs, i.e.
$S^{ai}=S^i$, $S^{bi}=(1-\la^i)S^i$,
where each $S^i$ is positive, $\la^i\in(0,1)$.
Denote the interval of the NGD prices in this model
by $I_\la(F)$ (the NGD pricing technique might be utility-based
or RAROC-based as the latter one is reduced to the former
one by considering
$\DDD=\frac{1}{1+R}\,\PD+\frac{R}{1+R}\,\RD$).
Let $(\la_n)=(\la_n;i\in I,n\in\N)$
be a sequence such that
$\la_n^i\xra[n\to\infty]{}0$ for any~$i$.

\begin{Theorem}
\label{PDT2}
For $F\in L_s^1(\DDD)$, we have
$I_{\la_n}(F)\xra[n\to\infty]{}I_0(F)$
in the sense that the right {\rm(}resp., left{\rm)}
endpoints of $I_{\la_n}(F)$ converge to the right
{\rm(}resp., left{\rm)} endpoint of $I_0(F)$.
\end{Theorem}

\sksminus
\Proof
Let $r$ denote the right endpoint of $I_0(F)$.
Suppose that the right endpoints of $I_{\la_n}(F)$ do not
converge to~$r$.
Then there exists $r'>r$ such that, for each~$n$
(possibly, after passing on to a subsequence),
there exists $\QQ_n\in\DDD\cap\RRR_{\la_n}$ with the property:
$\EE_\QQ F\ge r'$ ($\RRR_\la$ is the set of
risk-neutral measures in the model corresponding to~$\la$).
The sequence $(\QQ_n)$ has a
weak limit point $\QQ_\infty\in\DDD$.
Fix $i\in I$, $u\le v\in[0,T]$,
and a positive bounded $\F_u$-measurable function~$H$.
For any $n$, we have
$\EE_{\QQ_n} H((1-\la_n^i)S_v^i-S_u^i)\le0$.
As $S_v^i\in L_s^1(\DDD)$, we have
$\sup_{\QQ\in\DDD}\EE_\QQ S_v^i<\infty$, and hence,
$\limsup_n\EE_{\QQ_n}H(S_v^i-S_u^i)\le0$.
As the map
$\DDD\ni\QQ\mapsto\EE_\QQ H(S_v^i-S_u^i)$
is weakly continuous, we get
$\EE_{\QQ_\infty}H(S_v^i-S_u^i)\le0$.
In a similar way, we prove that
$\EE_{\QQ_\infty}H(-S_v^i+S_u^i)\le0$.
Thus, $S^i$ is an $(\F_t,\QQ_\infty)$-martingale, so that
$\QQ_\infty\in\DDD\cap\RRR_0$.
As the map $\DDD\ni\QQ\mapsto\EE_\QQ F$ is
weakly continuous, we should have
$\EE_{\QQ_\infty}F\ge r'$.
But this is a contradiction.\End

\subsection{Hedging}
\label{HI}

Consider the model of Subsection~\ref{UP}.

\begin{Definition}\rm
\label{HI1}
The \textit{upper and lower NGD prices} of a contingent
claim~$F$ are defined by
\begin{align*}
\overline V(F)&=\inf\{x:\exists X\in A\text{ such that }
u(X-F+x)\ge0\},\\
\underline V(F)&=\sup\{x:\exists X\in A\text{ such that }
u(X+F-x)\ge0\}.
\end{align*}
\end{Definition}

The problem of finding $\overline{V}(F)$
has some similarities with the
superreplication problem considered by Cvitani\'c,
Karatzas~\cite{CK99} and by Sekine~\cite{Sek04}, but the
difference is that in those papers the risk is measured
not as $\rho(X-F+x)$, but rather as $\rho((X-F+x)^-)$.

\begin{Proposition}
\label{HI2}
If $A$ is a cone and $F\in L_s^1(\DDD)$, then
\begin{align*}
\overline{V}(F)&=\sup\{\EE_\QQ F:\QQ\in\DDD\cap\RRR\},\\
\underline{V}(F)&=\inf\{\EE_\QQ F:\QQ\in\DDD\cap\RRR\}.
\end{align*}
\end{Proposition}

\Proof
Take $x_0\in\R$ and set
$A(x_0)=A+\{h(x_0-F):h\in\R_+\}$.
Using Theorem~\ref{UP4}, we can write
\begin{align*}
\wl V(F)\ge x_0
&\;\Lea\;\not\!\exists X\in A\text{ such that }u(X-F+x_0)>0\\
&\;\Lea\;\not\!\exists X\in A(x_0)\text{ such that }u(X)>0\\
&\;\Lea\;\DDD\cap\RRR(A(x_0))\ne\emp\\
&\;\Lea\;\exists\QQ\in\DDD\cap\RRR
\text{ such that }\EE_\QQ F\ge x_0.
\end{align*}
This yields the formula for $\overline{V}(F)$.
The representation of $\underline{V}(F)$ is proved
similarly.\End

\skm
\Remarks
(i) The above theorem is formally true if the NGD is violated.
In this case $\overline{V}(F)=-\infty$ and
$\underline{V}(F)=\infty$.

(ii) The above argument shows that there exist
$\overline{\QQ},\underline{\QQ}\in\DDD\cap\RRR$ such that
$\EE_{\overline{\QQ}}F=\overline{V}(F)$,
$\EE_{\underline{\QQ}}(F)=\underline{V}(F)$.
This is in contrast with the NA technique.\End

(iii) Under the conditions of the above corollary,
we have $\INGD(F)=[\underline{V}(F),\overline{V}(F)]$.

\skm
Let us now study the sub- and super-replication problem
for a particular case of a (frictionless) static model
with a finite number of assets.
Thus, we are given $S_0\in\R^d$ and
$S_1^1,\dots,S_1^d\in L_w^1(\DDD)$. From the financial
point of view, $S_n^i$ is the discounted price of the
$i$-th asset at time~$n$.

\begin{Definition}\rm
\label{HI3}
The \textit{superhedging and subhedging strategies}
are defined by
\begin{align*}
\overline H(F)&=\{h\in\R^d:
u(\lb h,S_1-S_0\rb-F+\overline V(F))\ge0\},\\
\underline H(F)&=\{h\in\R^d:
u(\lb h,S_1-S_0\rb+F-\underline V(F))\ge0\}.
\end{align*}
\end{Definition}

Below we provide a simple geometric procedure to
determine these quantities.
Assume that $F\in L_w^1(\DDD)$ and let us introduce the
notation
\begin{align*}
G&=\cl\{\EE_\QQ(S_1,F):\QQ\in\DDD\},\\
\overline v&=\sup\{x:(S_0,x)\in G\},\\
\underline v&=\inf\{x:(S_0,x)\in G\},\\
\overline N&=\{h\in\R^{d+1}:\forall x\in G,\,
\lb h,x-(S_0,\overline v)\rb\ge0\},\\
\underline N&=\{h\in\R^{d+1}:\forall x\in G,\,
\lb h,x-(S_0,\underline v)\rb\ge0\},
\end{align*}
i.e. $G$ is the generator for $(S_1,F)$ and $u$;
$\overline N$ (resp., $\underline N$) is the set
of inner normals to $G$ at the point $(S_0,\overline v)$
(resp., $(S_0,\underline v)$); see Figure~4.

\begin{Proposition}
\label{HI4}
We have
\begin{align*}
&\overline V(F)=\overline v,\\
&\underline V(F)=\underline v,\\
&\overline H=\{h\in\overline N:h^{d+1}=-1\},\\
&\underline H=\{h\in\underline N:h^{d+1}=1\}.
\end{align*}
\end{Proposition}

\Remark
The statement is true both in the case, where the NGD is
satisfied, and in the case, where the NGD is not satisfied
(in the latter case
$S_0$ does not belong to the projection of $G$ on $\R^d$,
$\overline v=\overline V(F)=-\infty$,
$\underline v=\underline V(F)=\infty$,
$\overline N=\overline H=\emp$,
$\underline N=\underline H=\emp$).

\skm
\texttt{Proof of Proposition~\ref{HI4}.}
This is an easy consequence of the line
$$
u(\lb h,S_1-S_0\rb\pm F\mp x)
=\inf_{z\in G}\lb(h,\pm1),z-(S_0,x)\rb,\quad h\in\R^d.
$$

\begin{figure}[h]
\begin{picture}(150,67.5)(-57.5,-20)
\put(5,-0.7){\includegraphics{poem.4}}
\put(0,0){\vector(1,0){45}}
\put(0,0){\vector(0,1){45}}
\multiput(20,-0.5)(0,2){23}{\line(0,1){1}}
\multiput(0,38.5)(1.9,0){11}{\line(1,0){1}}
\multiput(0,8.7)(1.9,0){11}{\line(1,0){1}}
\multiput(15.7,0)(0,2){10}{\line(0,1){1}}
\multiput(25.5,0)(0,2){15}{\line(0,1){1}}
\multiput(15.8,18.7)(1.7,0){3}{\line(1,0){0.9}}
\multiput(20,28.7)(2,0){3}{\line(1,0){1}}
\put(17.3,32.5){$\biggl\{$}
\put(19.2,12.7){$\biggr\}$}
\put(16,32.5){\small $1$}
\put(21.8,12.7){\small $1$}
\put(43,-4){\small $\R^d$}
\put(-4,43){\small $\R$}
\put(-9,8.5){\scalebox{0.8}{$\underline{V}(F)$}}
\put(-9,37.5){\scalebox{0.8}{$\overline{V}(F)$}}
\put(7,1.5){\scalebox{0.8}{$\underline{H}(F)$}}
\put(26,1.5){\scalebox{0.8}{$\overline{H}(F)$}}
\put(18.5,-4){\scalebox{0.8}{$S_0$}}
\put(30,33){\small $G$}
\put(-3.5,-15){\parbox{52mm}{\small\textbf{Figure~4.} Solution of the
super- and subhedging problem}}
\end{picture}
\end{figure}

\begin{Example}\rm
\label{HI5}
Consider the setting of Example~\ref{PF1}.
The results of that example show that
\begin{align*}
\overline V(F)&=\lb b,S_0-a\rb+\EE F+\al,\\
\underline V(F)&=\lb b,S_0-a\rb+\EE F-\al,
\end{align*}
where $\al=(\si^2\ga^2-\si^2\lb S_0-a,C^{-1}(S_0-a)\rb)^{1/2}$.
In order to find $\overline H$ and $\underline H$,
we express the upper and lower borders of the set~$\wt G$
given by~\eqref{pf0} as
$y=\pm(\si^2\ga^2-\si^2\lb x-a,C^{-1}(x-a)\rb)^{1/2}$.
Then by differentiation we get
$$
\overline H(\wt F)
=\underline H(\wt F)
=d|_{x=S_0}(\si^2\ga^2-\si^2\lb x-a,C^{-1}(x-a)\rb)^{1/2}
=-\si^2\al^{-1}C^{-1}(S_0-a).
$$
Hence,
\begin{align*}
\overline H(F)&=b-\si^2\al^{-1}C^{-1}(S_0-a),\\
\underline H(F)&=-b-\si^2\al^{-1}C^{-1}(S_0-a).
\end{align*}
\end{Example}

We will now provide one more example.
Let $S_0\in(0,\infty)$ and $S_1$ be an integrable random
variable such that $\Law S_1$ has no atoms and
$\supp\Law S_1=\R_+$.
Let $u$ be the coherent utility function corresponding
to Tail V@R of order $\la\in(0,1]$
(see Example~\ref{BD5}~(i)).
We assume that $u(S_1)<S_0<-u(-S_1)$.
Finally, let $F=f(S_1)$, where $f:\R_+\to\R$ is a
convex function of linear growth.
Let us denote $\Law S_1$ by $\QQ$ and find
$a$, $b$, $c$, $d$ such that
$a+b=\la$, $d-c=\la$, and
\begin{align*}
&\la^{-1}\int_0^{q_a}x\QQ(dx)
+\la^{-1}\int_{q_{1-b}}^\infty x\QQ(dx)=S_0,\\
&\la^{-1}\int_{q_c}^{q_d}x\QQ(dx)=S_0,
\end{align*}
where $q_x$ is the $x$-quantile of $\QQ$.

\begin{Proposition}
\label{EX1}
We have
\begin{align*}
\overline V(F)&=\la^{-1}\int_0^{q_a}f(x)\QQ(dx)
+\la^{-1}\int_{q_{1-b}}^\infty f(x)\QQ(dx),\\[1mm]
\underline V(F)&=\la^{-1}\int_{q_c}^{q_d}f(x)\QQ(dx),\\[1mm]
\overline H(F)&=\frac{f(q_{1-b})-f(q_a)}{q_{1-b}-q_a},\\[1mm]
\underline H(F)&=-\frac{f(q_d)-f(q_c)}{q_d-q_c}.
\end{align*}
\end{Proposition}

\Proof
Let us first prove the representation for $\overline V(F)$
under an additional assumption that $f$ is strictly
convex.
By Proposition~\ref{HI4},
$$
\overline V(F)=\sup_{Z\in\DDD_\la:\,\EE ZX=S_0}\EE Z f(X)
$$
($\DDD_\la$ is given by~\eqref{bd3}). Take
$$
Z_0\in\argmax_{Z\in\DDD_\la:\,\EE ZX=S_0}\EE Z f(X)
$$
($Z_0$ exists by a compactness argument).
Passing from $Z_0$ to $\EE(Z_0\mid X)$, we can assume
that $Z_0$ is $X$ measurable, i.e. $Z_0=\phi(X)$.
Let us prove that
\begin{equation}
\label{ex1}
Z_0=\la^{-1}I(X<q_a)+\la^{-1}I(X>q_{1-b}).
\end{equation}
Assume the contrary.
Then there exist $0<\al_1<\al_2<\al_3<\al_4$ such that
\begin{align*}
&\QQ(\{\phi<\la^{-1}\}\cap(\al_1,\al_2))>0,\\
&\QQ(\{\phi>0\}\cap(\al_2,\al_3))>0,\;\;\\
&\QQ(\{\phi<\la^{-1}\}\cap(\al_3,\al_4))>0.
\end{align*}
For $h_1,h_2,h_3\in[0,\la^{-1}]$, we set
$$
\wt\phi(x)=\begin{cases}
\phi(x),&x\notin(\al_1,\al_4),\\
\phi(x)\vee h_1,&x\in(\al_1,\al_2),\\
\phi(x)\wedge h_2,&x\in(\al_2,\al_3),\\
\phi(x)\vee h_3,&x\in(\al_3,\al_4).
\end{cases}
$$
We can find $h_1,h_2,h_3$ such that
\begin{align*}
&\hspace*{15mm}\QQ(\{\wt\phi>\phi\}\cap(\al_1,\al_2))>0,\\
&\hspace*{15mm}\QQ(\{\wt\phi<\phi\}\cap(\al_2,\al_3))>0,\\
&\hspace*{15mm}\QQ(\{\wt\phi>\phi\}\cap(\al_3,\al_4))>0,\\
&\int_0^\infty x\wt\phi(x)\QQ(dx)
=\int_0^\infty x\phi(x)\QQ(dx)=S_0,\\
&\int_0^\infty\wt\phi(x)\QQ(dx)
=\int_0^\infty\phi(x)\QQ(dx)=1.
\end{align*}
Consider the affine function $\wt f$ that coincides with
$f$ at $\al_2$ and $\al_3$. Then
$$
\int_0^\infty(\wt\phi(x)-\phi(x))\wt f(x)\QQ(dx)=0.
$$
Furthermore, as $f$ is strictly convex,
$\wt f<f$ on $(\al_1,\al_2)$,
$\wt f>f$ on $(\al_2,\al_3)$,
and $\wt f>f$ on $(\al_3,\al_4)$.
Consequently,
$$
\int_0^\infty(\wt\phi(x)-\phi(x))f(x)\QQ(dx)>0.
$$
Thus, we have found $\wt Z_0=\wt\phi(X)\in\DDD_\la$ such
that $\EE\wt Z_0 X=S_0$ and
$\EE\wt Z_0 f(X)>\EE Z_0 f(X)$, which
contradicts the choice of~$Z_0$.
As a result, \eqref{ex1} is satisfied, which yields the
desired representation of $\overline V(F)$.

Let us now prove the representation for $\overline{V}(F)$
in the general case.
Take $Z_0$ given by~\eqref{ex1}.
Find a strictly convex function $\wt f$ of linear growth.
Then the function $f_\eps=f+\eps\wt f$ is strictly convex
and the result proved above shows that
$\EE Zf_\eps(X)\le\EE Z_0f_\eps(X)$
for any $Z\in\DDD_\la$.
Passing on to the limit as $\eps\da0$,
we get $\EE Zf(X)\le\EE Z_0f(X)$
for any $Z\in\DDD_\la$.
This yields the desired representation of~$\overline V(F)$.

Let us now prove the representation for $\overline H(F)$.
Consider the function
$$
g(x)=\sup_{Z\in\DDD_\la:\,\EE ZX=x}\EE Zf(X),
\quad x\in[u(S_1),-u(-S_1)].
$$
It follows from the reasoning given above that
$g=g_1\circ g_2^{-1}$, where
\begin{align*}
g_1(x)&=\la^{-1}\int_0^{q_x}f(y)\QQ(dy)
+\la^{-1}\int_{q_{1-\la+x}}^\infty f(y)\QQ(dy),\quad
x\in[0,\la^{-1}],\\
g_2(x)&=\la^{-1}\int_0^{q_x}y\QQ(dy)
+\la^{-1}\int_{q_{1-\la+x}}^\infty y\QQ(dy),\quad
x\in[0,\la^{-1}].
\end{align*}
Applying Proposition~\ref{HI4}, we get
$$
\overline{H}(F)
=g'(S_0)
\frac{f(q_{1-b})-f(q_a)}{q_{1-b}-q_a}.
$$

The representations for $\underline V(F)$ and
$\underline H(F)$ are proved in a similar way.\End


\end{document}